\documentclass{agtart_a}
\pdfoutput=1

\usepackage{pinlabel}

\title{Sutured Heegaard diagrams for knots}
\author{Yi Ni}
\givenname{Yi}
\surname{Ni}
\address{Department of Mathematics\\Princeton University\\
Princeton\\NJ 08544\\USA}
\email{yni@math.princeton.edu}
\urladdr{}

\volumenumber{6}
\issuenumber{}
\publicationyear{2006}
\papernumber{18}
\startpage{513}
\endpage{537}

\doi{}
\MR{}
\Zbl{}

\keyword{knot Floer homology}
\keyword{sutured Heegaard diagram}
\keyword{Murasugi sum}
\keyword{semifibred satellite knot}
\subject{primary}{msc2000}{57R58}
\subject{primary}{msc2000}{57M27}
\subject{secondary}{msc2000}{53D40}

\received{22 July 2005}
\revised{13 January 2006}
\accepted{12 March 2006}
\published{2 April 2006}
\publishedonline{2 April 2006}
\proposed{}
\seconded{}
\corresponding{}
\editor{}
\version{}

\arxivreference{math.GT/0507440}

\AtBeginDocument{\let\bar\wbar\let\tilde\wtilde\let\hat\what\let\widetilde\wtilde\let\overline\wbar}

\def\int{\mathrm{int}}

\newcommand{\homo}{{\mathrm H}}

\makeatletter
\def\cnewtheorem#1[#2]#3{\newtheorem{#1}{#3}[section]
\expandafter\let\csname c@#1\endcsname\c@thm}

\newtheorem{thm}{Theorem}[section]
\cnewtheorem{lem}[thm]{Lemma}
\cnewtheorem{cor}[thm]{Corollary}
\cnewtheorem{prop}[thm]{Proposition}

\theoremstyle{definition}
\cnewtheorem{defn}[thm]{Definition}
\cnewtheorem{rem}[thm]{Remark}
\cnewtheorem{construction}[thm]{Construction}

\makeautorefname{defn}{Definition}

\makeatother  

\begin{document}

\begin{asciiabstract}
We define sutured Heegaard diagrams for null-homologous knots in
3-manifolds. These diagrams are useful for computing the knot Floer
homology at the top filtration level.  As an application, we give a
formula for the knot Floer homology of a Murasugi sum. Our result
echoes Gabai's earlier works. We also show that for so-called
`semifibred' satellite knots, the top filtration term of the knot
Floer homology is isomorphic to the counterpart of the companion.
\end{asciiabstract}

\begin{abstract}
We define sutured Heegaard diagrams for null-homologous knots in
3--manifolds. These diagrams are useful for computing the knot Floer
homology at the top filtration level. As an application, we give a
formula for the knot Floer homology of a Murasugi sum. Our result echoes
Gabai's earlier works. We also show that for so-called
``semifibred" satellite knots, the top filtration term of the knot
Floer homology is isomorphic to the counterpart of the companion.
\end{abstract}

\maketitle

\section{Introduction}

Knot Floer homology was introduced by Ozsv\'ath and Szab\'o in
\cite{OSz2}, and independently by Rasmussen in \cite{Ra}, as part of
Ozsv\'ath and Szab\'o's Heegaard Floer theory. A survey of Heegaard
Floer theory can be found in Ozsv\'ath and Szab\'o \cite{OSz9}.

One remarkable feature of knot Floer homology is that it determines
the genus in the case of classical knots (Ozsv\'ath and Szab\'o
\cite[Theorem 1.2]{OSz4}), namely, the genus of a classical knot is
the highest nontrivial filtration level of the knot Floer
homology. The proof of this deep result uses Gabai's work on the
existence of taut foliations of knot complements \cite{G3}. Another
theorem of Gabai can be used to generalize Ozsv\'ath and Szab\'o's
result to links in homology 3--spheres (Ni \cite{Ni}). Hence one may
naturally expect that there is a more precise relationship between
taut foliation and the top filtration term of knot Floer homology.

Another interesting property of knot Floer homology is that, for
fibred knots, the top filtration term of knot Floer homology is a
single $\mathbb Z$ (Ozsv\'ath and Szab\'o \cite[Theorem
1.1]{OSz3}). It is conjectured that the converse is also true for
classical knots.

The results cited above show that a lot of information about the
knot is contained in the top filtration term of knot Floer
homology.  In the present paper, we introduce sutured Heegaard
diagrams for knots, which are useful for computing the top
filtration term of knot Floer homology.

The definition of a sutured Heegaard diagram will be given in Section
2.  We state here two theorems as applications:

\begin{thm}
Suppose $K_1,K_2\subset S^3$ are two knots and that $K$ is the
Murasugi sum performed along the minimal genus Seifert
surfaces. Suppose the genera of $K_1,K_2,K$ are $g_1,g_2,g$,
respectively, then
$$\widehat{HFK}(K,g;\mathbb F)\cong\widehat{HFK}(K_1,g_1;\mathbb F)\otimes\widehat{HFK}(K_2,g_2;\mathbb F)$$
as linear spaces, for any field $\mathbb F$.
\end{thm}

\begin{thm}
Suppose that $K^*$ is a semifibred satellite knot of $K$. Suppose the
genera of $K,K^*$ are $g,g^*$, respectively. Then
$$\widehat{HFK}(K^*,g^*)\cong\widehat{HFK}(K,g)$$
as abelian groups.
\end{thm}

The precise definitions of Murasugi sum and semifibred satellite
knot will be given later, where we will prove more general
versions of these two theorems.

The paper is organized as follows.

In Section 2 we review the adjunction inequality. Then we give the
definition of a sutured Heegaard diagram.

In Section 3, we enhance Ozsv\'ath and Szab\'o's winding argument.
Using this argument, we show that a sutured Heegaard diagram may
conveniently be used to compute the top filtration term of the knot
Floer homology. As an immediate application, we give a new proof of a
result due to Ozsv\'ath and Szab\'o.

Section 4 will be devoted to the study of Murasugi sum. The
formula for Murasugi sum is almost a direct corollary of the
results in Section 3, once we know what the Heegaard diagram is.
Our formula echoes Gabai's earlier works.

Section 5 is about semifibred satellite knots. Again, most efforts
are put on the construction of a Heegaard diagram.

\medskip
\noindent{\bf Acknowledgements}\qua We are grateful to David Gabai,
Jacob Rasmussen and Zolt\'an Szab\'o for many stimulating discussions
and encouragements. We are especially grateful to the referee for a
detailed list of corrections and suggestions.

The author is partially supported by the Centennial fellowship of the
Graduate School at Princeton University. Part of this work was carried
out during a visit to Peking University; the author wishes to thank
Shicheng Wang for his hospitality during the visit.

\section{Sutured Heegaard diagrams}

The definition of a sutured Heegaard diagram relies on a detailed
understanding of the adjunction inequality for knot Floer homology.
For the reader's convenience, we sketch a proof here, which is derived
from arguments of Ozsv\'ath and Szab\'o (\cite[Theorem 7.1]{OSz1} and
\cite[Theorem 5.1]{OSz2}).

\begin{nopartheorem}
\begin{thm}[\protect{\cite[Theorem 7.1]{OSz1}}]\label{adjunct}
Let $K\subset Y$ be an oriented null-homologous knot, and suppose
that $\widehat{HFK}(Y,K,\underline{\mathfrak{s}})\ne0$. Then, for
each Seifert surface F for K of genus $g>0$, we have that
$$\left|\langle c_1(\underline{\mathfrak{s}}),[\widehat{F}]\rangle\right|\le2g(F).$$
\end{thm}\end{nopartheorem}

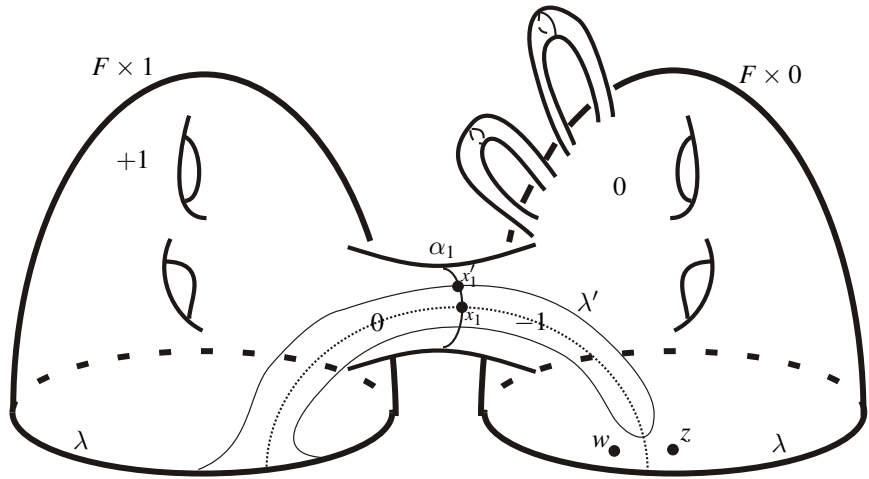
\begin{figure}[ht!]\small\label{Fig:1}
\labellist
\pinlabel $F\times1$ [br] at 138 787
\pinlabel $F\times0$ [bl] at 410 783
\pinlabel $0$ at 356 738
\pinlabel $0$ at 241 673
\pinlabel $-1$ at 314 673
\pinlabel $\lambda$ [bl] at 94 610
\pinlabel $\lambda$ [br] at 438 609
\pinlabel $\lambda'$ [bl] at 333 676
\pinlabel $\alpha_1$ [b] at 272 700
\pinlabel* {$\scriptstyle x_1$} [tl] <1.5pt,-1.5pt> at 281 680
\pinlabel* {$\scriptstyle x_1'$} [bl] <1.5pt,0pt> at 279 690
\pinlabel $w$ [br] at 354 611
\pinlabel $+1$ at 125 750
\pinlabel $z$ [bl] <1pt,1pt> at 381 612
\endlabellist
\cl{\includegraphics[width=4.5in]{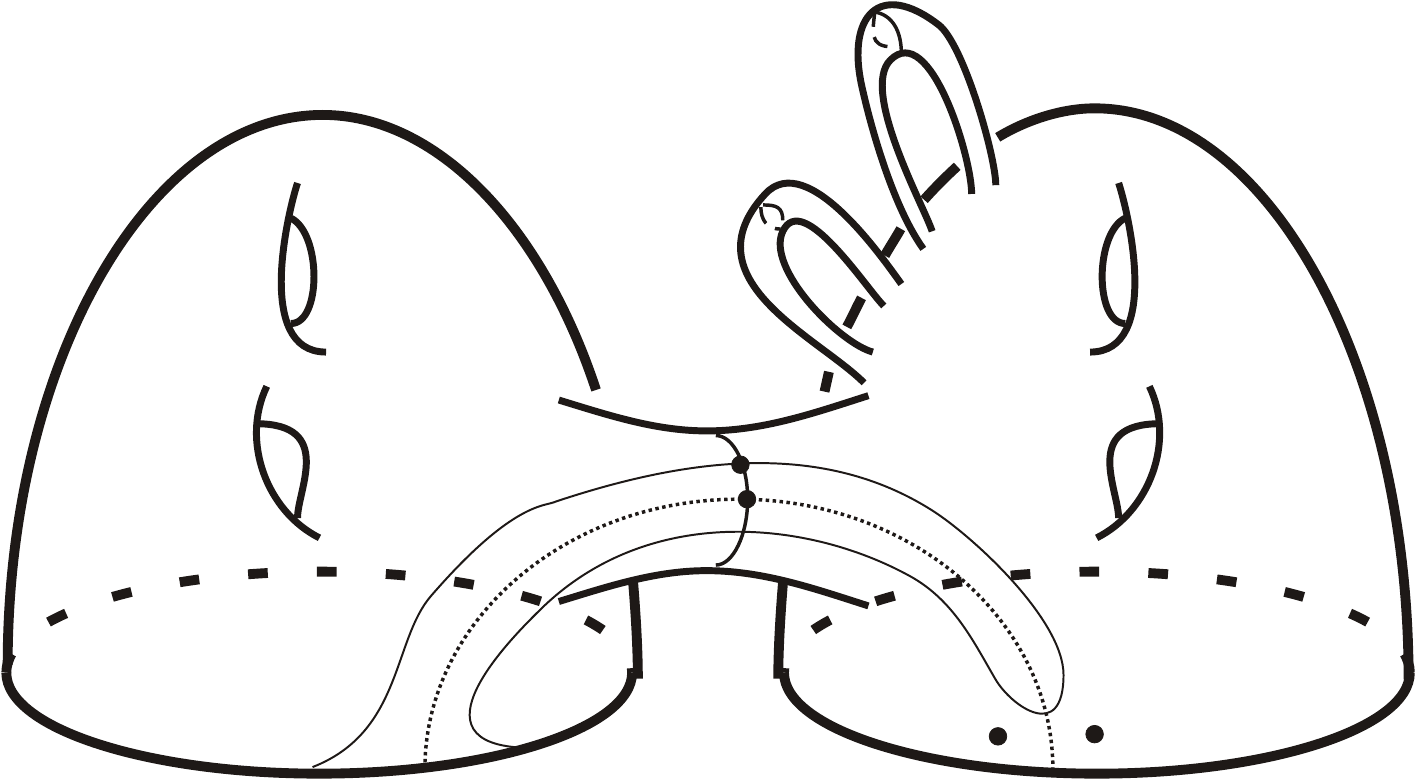}}
\caption{The Heegaard surface $\Sigma$, here the two $\lambda$'s are
glued together. Note the multiplicities of $\mathcal P$ in the
regions on the diagram, denoted by $-1,0,+1$.}
\end{figure}

\begin{proof}[Sketch of proof] The proof consists of 3 steps.

{\bf Step 1}\qua {\sl Construct a Heegaard splitting of $Y$}

Consider a product neighborhood of $F$ in $Y$:
$N(F)=F\times[0,1]$. Let $${\mathcal H}=D^1\times D^2\subset
Y-\int(N(F))$$ be a 1--handle connecting $F\times0$ to $F\times1$,
$\partial_v{\mathcal H }=D^1\times\partial D^2$ is the vertical
boundary of ${\mathcal H}$. We can choose ${\mathcal H}$ so that
it is ``parallel" to $\textsl{point}\times[0,1]$ for a point on
$\partial F$. Namely, there is a properly embedded product disk
${\mathfrak D}$ in $Y-\int(N(F)\cup{\mathcal H})$,
$\partial{\mathfrak D}\cap(\partial F\times[0,1])$ is an essential
arc in $\partial F\times[0,1]$, and $\partial{\mathfrak D
}\cap\partial_v{\mathcal H}$ is an essential arc in
$\partial_v{\mathcal H}$.

Let $M=Y-\int(N(F)\cup{\mathcal H}\cup N(\mathfrak D))$, $f\co M\to
[0,3]$ be a Morse function, so that $\partial M=f^{-1}(0)$,
$\{\textrm{Critical points of index } i\}\subset f^{-1}(i)$ for
$i=0,1,2,3$. Let $\widetilde{\Sigma}=f^{-1}(\frac32)$,
$\widetilde{\Sigma}$ separates $Y$ into two parts $\widetilde
U_0,\widetilde U_1$. Suppose $\widetilde U_0$ is the part
containing $F$, then $\widetilde U_0$ can be obtained by adding
$r$ 1--handles to $N(F)\cup{\mathcal H}\cup N(\mathfrak D)$. After
handlesliding, one can assume these 1--handles are attached to
$F\times0$. Moreover, we can let the two feet of any 1--handle be
very close on $F\times0$.

$\widetilde U_1$ is a handlebody, and $N(\mathfrak D)$ is a
1--handle attached to it, if you turn $\widetilde
U_0\cup\widetilde U_1$ upsidedown. Now let $U_1=\widetilde U_1\cup
N(\mathfrak D)$, $U_0=Y-\int(U_1)$, then $Y=U_0\cup_{\Sigma}U_1$ is
a Heegaard splitting of $Y$.

\medskip
{\bf Step 2}\qua {\sl Find a weakly admissible Heegaard diagram for
$(Y,K)$}

For each 1--handle $D^1\times D^2$ attached to $F\times[0,1]$,
choose its belt circle $0\times \partial D^2$ ($D^1$ viewed as
$[-1,1]$). $\alpha_1$ denotes the belt circle of $\mathcal H$, and
the belt circles of other 1--handles are denoted by
$\alpha_{2g+2},\dots,\alpha_{2g+1+r}$.

Choose a set of disjoint, properly embedded arcs $\xi_i$
($i=2,3,\dots,2g+1$) on $F\times1$, so that they represent a basis
of $\homo_1(F,\partial F)$. Choose a copy of $\xi_i$ on
$F\times0$, denoted by $\overline{\xi_i}$. For each $i$, complete
$\xi_i\sqcup\overline{\xi_i}$ by two vertical arcs on $\partial
F\times[0,1]$ to get a simple closed curve $\alpha_i$ on
$\partial(F\times[0,1])$, which bounds a disk in $F\times[0,1]$.
($\alpha_i$ can be viewed as the ``double" of $\xi_i$.) We can
choose $\xi_i$ so that $\alpha_i$ is disjoint from $\partial
\mathfrak D$, $i=2,\dots,2g+1$.

Let $\lambda=\partial F\times1$ be the longitude of $K$, and
$\mu=\partial\mathfrak D$ will be the meridian of $K$. Both
$\lambda$ and $\mu$ are simple closed curves on $\Sigma$. Extend
$\mu$ to a set of disjoint simple closed curves
$\{\mu,\beta_2,\beta_3,\dots,\beta_{2g+1+r}\}$ on $\Sigma$, so
that they are linearly independent in $\homo_1(\Sigma)$, and each
bounds a non-separating disk in $U_1$. Let
$\mbox{\boldmath${\alpha}$}=\{\alpha_1,\dots,\alpha_{2g+1+r}\},\mbox{\boldmath$\beta_0$}=\{\beta_2,\dots,\beta_{2g+1+r}\}$.
 Then
$(\Sigma,\mbox{\boldmath${\alpha}$},\mbox{\boldmath$\beta_0$}\cup\{\mu\})$
is a Heegaard diagram for $Y$.

Since $|\lambda\cap\mu|=1$, by handleslides over $\mu$, we can
assume all $\beta_i$'s ($2\le i\le 2g+1+r$) are disjoint from
$\lambda$. Now it is easy to see
$(\Sigma,\mbox{\boldmath${\alpha}$},
\mbox{\boldmath$\beta_0$},\mu,\lambda\cap\mu)$ is a marked
Heegaard diagram for the knot $(Y,K)$. Furthermore, we construct a
double-pointed Heegaard diagram
$(\Sigma,\mbox{\boldmath${\alpha}$},
\mbox{\boldmath$\beta_0$}\cup\{\mu\},w,z)$ from the marked one.

Choose a set of circles $\{\tau_2,\tau_3,\dots,\tau_{2g+1+r}\}$ on
$\Sigma-\lambda-\alpha_1-\mu$, so that $\tau_i$ intersects
$\alpha_i$ transversely and exactly once,
$\tau_i\cap\alpha_j=\emptyset$ when $i\ne j$. Wind $\alpha_i$'s
along $\tau_i$'s ($i=2,\dots,2g+1+r$), one can get a weakly
admissible Heegaard diagram for $(Y,K)$. For more details, see
\cite[Theorem 7.1]{OSz1} or the discussion after \fullref{relative}.

\medskip
{\bf Step 3}\qua {\sl Proof of the inequality}\qua

On $\Sigma$, there is a domain $\mathcal P$ bounded by $\alpha_1$
and $\lambda$, which is basically $F\times1$ with a hole. Move
$w,z$ slightly out of $\mathcal P$. Hence $\mathcal P$ is a
periodic domain for the $Y_0$ Heegaard diagram
$(\Sigma,\mbox{\boldmath${\alpha}$},
\mbox{\boldmath$\beta_0$}\cup\{\lambda\},w)$.

Wind $\lambda$ once along $\mu$ as shown in \fullref{Fig:1}, which
should be compared with \cite[Figure~6]{OSz2}.  After winding,
$\lambda$ becomes a new curve $\lambda'$, and $\mathcal P$ becomes
$\mathcal P'$. By our choice, $\mu$ has no intersection with
$\xi_i$. Hence any intersection point $\mathbf x$ for the $(Y,K)$
diagram $(\Sigma,\mbox{\boldmath${\alpha}$},
\mbox{\boldmath$\beta_0$}\cup\{\mu\},w,z)$ must contain
$x_1=\mu\cap\alpha_1$. Let $\mathbf x'$ be the nearby intersection for
the $Y_0$ diagram $(\Sigma,\mbox{\boldmath${\alpha}$},
\mbox{\boldmath$\beta_0$}\cup\{\lambda'\},w)$, $\mathbf x'$ contains
$x_1'$. Here $Y_0$ is the manifold obtained from $Y$ by 0--surgery on
$K$. Local multiplicity of $\mathcal P'$ at $x_1'$ is 0. The same
argument as in the proof of \cite[Theorem 5.1]{OSz2} shows that
$$\langle c_1(\mathfrak{s}'({\mathbf x}')),[\widehat F]\rangle=-2g+\#(x_i \ \textrm{in the interior of}\ {\mathcal P})\ge-2g.$$
By conjugation invariance, we have the adjunction inequality.
\end{proof}

The reader should note, in the proof of 
\cite[Theorem 5.1]{OSz2}, it is not assumed that $\alpha_1$ is the only
$\alpha$ curve intersecting $\mu$.

Before moving on, we clarify one convention we are going to use.
The boundary map in Heegaard Floer theory is defined via counting
holomorphic disks in $\textrm{Sym}^n\Sigma$. There is a natural
$n$--value map $\varrho\co\textrm{Sym}^n\Sigma\to\Sigma$. In
practice we always consider the image of the holomorphic disk
under $\varrho$. By abuse of notation, we do not distinguish a
holomorphic disk and its image under $\varrho$.

As we have seen in the above proof, the chain complex
$\widehat{CFK}(Y,K,-g)$ is generated by the intersection points
$$\{\mathbf x\;|\;\mbox{no $x_i$ lies in the interior of $\mathcal
P$}\}$$ for the diagram $(\Sigma,\mbox{\boldmath${\alpha}$},
\mbox{\boldmath$\beta_0$}\cup\{\mu\},w,z)$. However, it is not
clear that the holomorphic disks connecting the generators of
$\widehat{CFK}(Y,K,-g)$ do not intersect the interior of $\mathcal
P$. In fact, if two components of
$\mbox{\boldmath$\beta_0$}\cap\mathcal P$ are parallel, then it is
very possible to have a holomorphic disk of quadrilateral type,
which connects two generators of $\widehat{CFK}(Y,K,-g)$ and
intersects the interior of $\mathcal P$.

However, if the Heegaard diagram is good enough, we can let all
the holomorphic disks be supported away from $\lambda$. This
observation leads to the following

\begin{defn}\label{suturedHD}
A double pointed Heegaard diagram
$$(\Sigma,\mbox{\boldmath${\alpha}$},
\mbox{\boldmath$\beta_0$}\cup\{\mu\},w,z)$$ for $(Y,K)$ is a {\it
sutured Heegaard diagram}, if it satisfies:

(Su0)\qua There exists a subsurface $\mathcal P\subset\Sigma$, bounded
by two curves $\alpha_1\in\mbox{\boldmath${\alpha}$}$ and
$\lambda$. $g$ denotes the genus of $\mathcal P$.

(Su1)\qua $\lambda$ is disjoint from \mbox{\boldmath$\beta_0$}. $\mu$
does not intersect any $\alpha$ curves except $\alpha_1$. $\mu$
intersects $\lambda$ transversely in exactly one point, and
intersects $\alpha_1$ transversely in exactly one point.
$w,z\in\lambda$ lie in a small neighborhood of $\lambda\cap\mu$,
and on different sides of $\mu$. (In practice, we often push $w,z$
off $\lambda$ into $\mathcal P$ or $\Sigma-\mathcal P$.)

(Su2)\qua $(\mbox{\boldmath$\alpha$}-\{\alpha_1\})\cap\mathcal P$
consists of $2g$ arcs, which are linearly independent in
$\homo_1(\mathcal P,\partial\mathcal P)$. Moreover,
$\Sigma-\mbox{\boldmath$\alpha$}-\mathcal P$ is connected.
\end{defn}

The existence of a sutured Heegaard diagram is guaranteed by the
construction in the proof of \fullref{adjunct}.

\section{Top filtration term of the knot Floer homology}

The following definition will be useful:

\begin{defn}\label{relative} In \cite[Section 2.4]{Osz0}, there are definitions of domain and
periodic domain. One can generalize these definitions to relative
case. For example, suppose $\mathcal R$ is an oriented connected
compact surface, with some base points on it. $\Gamma$ is a set of
finitely many properly embedded, mutually transverse curves on
$\mathcal R$. Let $\mathcal D_1,\dots,\mathcal D_m$ be the
closures of $\mathcal R-\Gamma$. A {\it relative domain} $\mathcal
D$ is a linear combination of the $\mathcal D_i$'s. A {\it
relative periodic domain} is a relative domain $\mathcal D$, such
that $\partial\mathcal D-\partial\mathcal R$ is a linear
combination of curves in $\Gamma$, and $\mathcal D$ avoids the
basepoints on $\mathcal R$. If a relative periodic domain has
nonnegative local multiplicities everywhere, then we call it a
{\it nonnegative} relative periodic domain.
\end{defn}

In order to get admissible diagrams, Ozsv\'ath and Szab\'o
introduced the technique of winding in \cite{Osz0}. We briefly
review this technique in our relative settings.

\begin{figure}[ht!]\small\label{Fig:2}
\labellist
\pinlabel $\tau$ [l] at 133 676
\pinlabel $a{-}1$ at 88 676
\pinlabel $a{-}1$ at 340 676
\pinlabel $a$ at 88 793
\pinlabel $a$ at 340 793
\pinlabel $\xi'$ [t] at 340 739
\pinlabel {Figure 2a} [t] <0pt, -4pt> at 133 642
\pinlabel {Figure 2b} [t] <0pt, -4pt> at 392 642
\pinlabel {Figure 2c} [t] <0pt, -4pt> at 133 369
\pinlabel {Figure 2d} [t] <0pt, -4pt> at 392 369
\endlabellist
\cl{\includegraphics[width=3.5in]{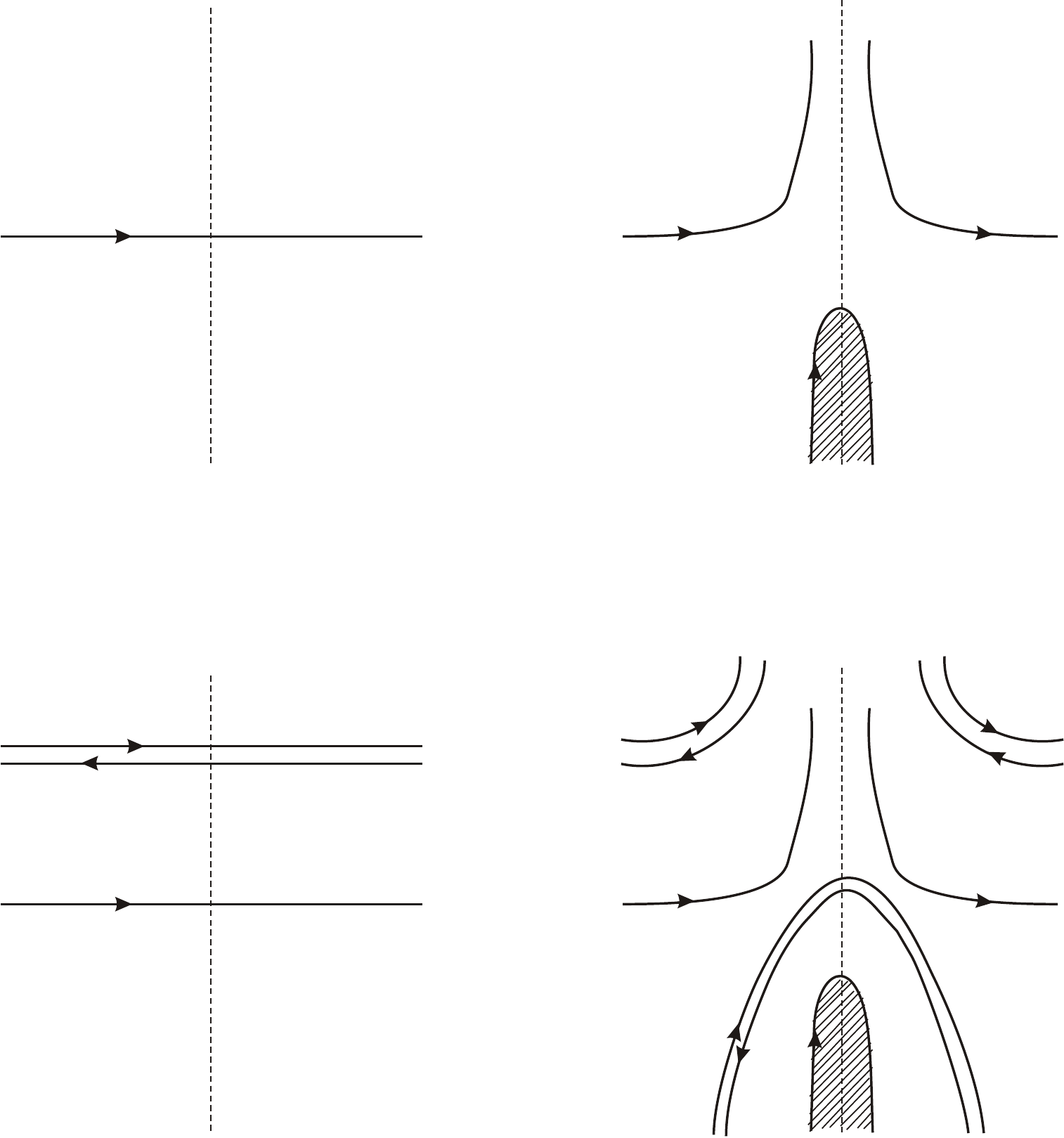}}\vspace{2mm}
\caption{Wind $\xi$ along $\tau$}
\end{figure}

Let $(\mathcal R,\Gamma)$ be as in \fullref{relative},
$\xi\in\Gamma$. Suppose there exists a simple closed curve
$\tau\subset\mathcal R$, which intersects $\xi$ transversely once.
We can wind $\xi$ once along $\tau$, as shown in \fullref{Fig:2}b. Now
suppose $\mathcal D$ is a relative periodic domain, such that
$\xi$ nontrivially contributes to $\partial\mathcal D$, say, the
contribution is 1.

In \fullref{Fig:2}a, the local multiplicities of $\mathcal D$ in the two
regions are $a$ and $a-1$, respectively. Suppose $\mathcal D'$ is
a variant of $\mathcal D$ after winding. In \fullref{Fig:2}b, the local
multiplicity of $\mathcal D'$ in the shaded area is $a-2$. If we
wind $\xi$ along $\tau$ sufficiently many times, we can get
negative local multiplicity here. If we wind $\xi$ along a
parallel copy of $\tau$, but in the other direction, we can get
positive local multiplicity.

Sometimes we have to wind several curves simultaneously along
$\tau$. We require the curve $\xi$, which we care about, has
nonzero algebraic intersection number with $\tau$; and all other
curves we wind have zero algebraic intersection number with
$\tau$. Then we can still get both positive and negative local
multiplicities after winding. See \fullref{Fig:2}c and \fullref{Fig:2}d for a
typical example.

Now we can give a key lemma:

\begin{lem}\label{winding}
Suppose $\mathcal R$ is an oriented connected compact surface,
with a nonempty collection of base points ${\mathbf w}$. Let
$\Gamma$ be a set of finitely many properly embedded, mutually
transverse curves on $\mathcal R$. $\Gamma=\Xi\sqcup \Theta$,
$\Xi=\{\xi_1,\dots,\xi_n\}$. Suppose that there are circles
$\tau_1,\dots,\tau_n$, such that $\tau_i$ intersects $\xi_i$
transversely in a single intersection point. Furthermore, suppose
there exists a region (connected open set) $\mathcal U$, so that
${\mathbf w}\subset\mathcal U$, $\tau_i\cap\xi_i\in\mathcal U$ for
all $i$, and all curves in $\Theta$ are disjoint from $\mathcal
U$.

Then after winding $\xi$ curves along the $\tau$ curves
sufficiently many times, every relative periodic domain $\mathcal
D$ whose boundary contains $\xi_i$ nontrivially (ie,
$n_i\cdot\xi_i\subset\partial\mathcal D$, $n_i\ne0$) has both
positive and negative local multiplicities. Hence there is {\bf
no} nonnegative relative period domain $\mathcal D$ with its
boundary containing $\xi_i$ nontrivially.

Furthermore, if the algebraic intersection number of $\tau_i$ with
$\xi_j$ is zero when $i\ne j$, and $\xi_i$'s are mutually
disjoint, then we can arrange that the $\xi_i$'s are mutually
disjoint after winding.
\end{lem}
\begin{proof}
For convenience, we use $\mathbb Q$ coefficients. Without loss of
generality, we can assume $\mathbf w$ consists of a single point
$w$. Let $\mathscr X$ be the linear space generated by curves in
$\Xi$, $\mathscr Y$ be the linear space generated by curves in
$\Theta$. There is a natural homomorphism
$$\mathscr H\co\mathscr X\oplus\mathscr Y\to\homo_1(\mathcal R,\partial\mathcal R).$$
For each nontrivial element $\gamma\in\ker\mathscr H$, there is a
unique relative periodic domain bounded by $\gamma$. Let $\mathscr
X_0=Proj_{\mathscr X}(\ker\mathscr H)$. Without loss of
generality, we can choose a basis of $\mathscr X_0$ in the form
$$\zeta_i=\xi_i+\sum_{j=m+1}^{n}c_{ij}\xi_j,\quad i=1,\dots,m.$$
Each $\zeta_i$ cobounds a relative periodic domain $\mathcal Q_i$
with some element in $\mathscr Y$.

We will wind $\xi_i$ along $\tau_i$ in one direction sufficiently
many times, and wind $\xi_i$ along a parallel copy of $\tau_i$ in
the other direction sufficiently many times, thus get new
collection of curves $\Xi'$. The variants of $\mathcal Q_i$ after
winding are denoted by $\mathcal Q_i'$. Then for each $i$, there
are points $w_i,z_i$ near $\xi_i\cap\tau_i$ (hence
$w_i,z_i\in\mathcal U$), such that $n_{w_i}(\mathcal Q_i')$ is
positive, $n_{z_i}(\mathcal Q_i')$ is negative.

Note that the winding along $\tau$ only changes the local
multiplicities in a neighborhood of $\tau$. We can choose those
neighborhoods narrow enough, so that $w_i,z_i$ are not in the
neighborhood of $\tau_j$ when $j\ne i$. Hence $n_{w_i}(\mathcal
Q_j')=n_{w_i}(\mathcal Q_j)$ when $j\ne i$. We wind $\xi_i$
sufficiently many times, so that
$$
n_{w_i}(\mathcal Q_i')>\sum_{j\ne i}|n_{w_i}(\mathcal Q_j)|,\qquad
|n_{z_i}(\mathcal Q_i')|>\sum_{j\ne i}|n_{z_i}(\mathcal Q_j)|.
$$
Now if $\mathcal D$ is a relative periodic domain, and
$\partial\mathcal D$ contains some $\xi_i'$ nontrivially, then
$Proj_{\mathscr X}(\partial\mathcal D)=\sum_{i=1}^{m}c_i\zeta_i'$
is a nontrivial sum. $\mathcal D -\sum_{i=1}^mc_i\mathcal Q_i'$ is
a relative periodic domain bounded by an element in $\mathscr Y$.
But all curves generating $\mathscr Y$ are disjoint from $\mathcal
U\ni w$, so $\mathcal D-\sum_{i=1}^mc_i\mathcal Q_i'$ is disjoint
from $\mathcal U$. Suppose $c_l$ is the coefficient with maximal
absolute value, then $\mathcal D$ has negative multiplicity at
$w_l$ or $z_l$.

If the algebraic intersection number of $\tau_i$ with $\xi_j$ is
zero when $i\ne j$, and $\xi_i$'s are mutually disjoint before
winding, then when we wind along $\tau_i$, we simultaneously wind
all the $\xi$ curves intersecting $\tau_i$. Hence $\xi_i$'s are
still disjoint after winding. We can get our result about local
multiplicity by the discussion before \fullref{winding}.
\end{proof}

Suppose $(\Sigma,\mbox{\boldmath${\alpha}$},
\mbox{\boldmath$\beta_0$}\cup\{\mu\},w,z)$ is a sutured Heegaard
diagram for $(Y,K)$. As in the proof of \fullref{adjunct}, the
generators of $\widehat{CFK}(Y,K,-g)$ are supported outside the
interior of $\mathcal P$. Our main result is

\begin{prop}\label{support}
Let $(\Sigma,\mbox{\boldmath${\alpha}$},
\mbox{\boldmath$\beta_0$}\cup\{\mu\},w,z)$ be a sutured Heegaard
diagram for $(Y,K)$. Then after winding transverse to the $\alpha$
curves, we get a new sutured Heegaard diagram
$$(\Sigma,\mbox{\boldmath${\alpha''}$},
\mbox{\boldmath$\beta_0$}\cup\{\mu\},w,z),$$ which is weakly
admissible, and all the holomorphic disks connecting generators of
$\widehat{CFK}(Y,K,-g)$ are supported outside a neighborhood of
$\lambda$.
\end{prop}
\begin{proof} Suppose the components of
$(\mbox{\boldmath$\alpha$}-\{\alpha_1\})\cap\mathcal P$ are
$\xi_2,\dots,\xi_{2g+1}$. Let $\mathcal U$ be a small neighborhood
of $\lambda$.

$\xi_2,\dots,\xi_{2g+1}$ are linearly independent in
$\homo_1(\mathcal P,\partial\mathcal P)$, and they are disjoint
from $\mu\cap\mathcal P$, which is an arc connecting $\alpha_1$ to
$\lambda$. Hence $\mathcal P-\cup_{i=2}^{2g+1}\xi_i-\mu$ is
connected. Now we can find simple closed curves
$\tau_2,\dots,\tau_{2g+1}\subset\mathcal P$, so that they are
disjoint from $\mu$ and $\xi$ curves, except that $\tau_i$
intersects $\xi_i$ transversely in a single intersection point. We
can assume $\tau_i\cap\xi_i\in\mathcal U$. By \fullref{winding},
after winding $\xi$ curves along $\tau$ curves, we can get a new
Heegaard diagram
$$(\Sigma,\mbox{\boldmath${\alpha'}$},
\mbox{\boldmath$\beta_0$}\cup\{\mu\},w,z),$$ such that its
restriction to $\mathcal P$ satisfies the conclusion of \fullref{winding}.

Suppose the closed $\alpha$ curves in $\Sigma-\mathcal P$ are
$\widetilde{\alpha}_{2g+2},\dots,\widetilde{\alpha}_{2g+1+r}$. By
Condition (Su2), $\Sigma-\mbox{\boldmath$\alpha$}-\mathcal P$ is
connected, hence we can find circles
$\tau_{2g+2},\dots,\tau_{2g+1+r}\subset(\Sigma-\mathcal P)$, such
that they are disjoint from all the $\alpha$ curves, except that
$\tau_i$ intersects $\widetilde{\alpha}_i$ transversely in a
single intersection point.

Now if $\mathcal D$ is a periodic domain, then $\mathcal
D\cap\mathcal P$ is a relative periodic domain in $\mathcal P$.
Hence if $\mathcal D$ is non-negative, then $\partial\mathcal D$
does not pass through $\xi_2,\dots,\xi_{2g+1}$. A similar argument
as in the proof of \fullref{winding} shows that we can wind
$\widetilde{\alpha}_{2g+2},\dots,\widetilde{\alpha}_{2g+1+r}$
along $\tau_{2g+2},\dots,\tau_{2g+1+r}$, to get a new diagram
$$(\Sigma,\mbox{\boldmath${\alpha''}$},
\mbox{\boldmath$\beta_0$}\cup\{\mu\},w,z),$$ such that $\partial
\mathcal D$ does not pass through
$\widetilde{\alpha}_{2g+2},\dots,\widetilde{\alpha}_{2g+1+r}$.

So $\partial\mathcal D$ consists of $\alpha_1$ and curves in
$\mbox{\boldmath$\beta_0$}$. Moreover, consider the curve $\mu$.
We observe that the base points $w,z$ lie close to and on both
sides of $\mu$, and $\mu$ has only one intersection with the
$\alpha$ curves. So $\mathcal D\cap\mu=\emptyset$. $\mu$
intersects $\alpha_1$, hence $\partial\mathcal D$ does not contain
$\alpha_1$. We conclude that
$(\Sigma,\mbox{\boldmath${\alpha''}$},
\mbox{\boldmath$\beta_0$}\cup\{\mu\},w,z)$ is weakly admissible,
since curves in $\mbox{\boldmath$\beta_0$}$ are linearly
independent in $\homo_1(\Sigma)$.

If $\Phi$ is a holomorphic disk connecting two generators of
$\widehat{CFK}(Y,K,-g)$, then $\Phi\cap\mathcal P$ is a relative
periodic domain in $\mathcal P$, since the generators lie outside
$\int(\mathcal P)$. So $\partial\Phi$ does not pass through $\xi$
curves. Hence $\Phi$ is disjoint from $\lambda$, since it should
avoid $w,z$.
\end{proof}

\begin{rem}\label{support2}
In practice, in order to compute $\widehat{HFK}(Y,K,-g)$ from a
given sutured Heegaard diagram, we only need to wind the closed
$\alpha$ curves in $\Sigma-\mathcal P$ sufficiently many times,
then count the holomorphic disks which are disjoint from
$\lambda$. The reason is that these disks are not different from
those disks obtained after winding $\xi$ curves.
\end{rem}

As an immediate application, we give the following proposition.
This proposition and the second proof here were told to the author
by Zolt\'an Szab\'o. The current paper was partially motivated by
an attempt to understand this proposition.

\begin{prop}\label{reglue}
$(Y,K)$ is a knot, $F$ is its Seifert surface of genus $g$. We cut
open $Y$ along $F$, then reglue by a self-diffeomorphism $\varphi$
of $F$. Denote the new knot in the new manifold by $(Y',K')$. Then
$$\widehat{HFK}(Y,K,-g)\cong\widehat{HFK}(Y',K',-g)$$
as abelian groups.
\end{prop}
\begin{proof}[The first proof] Construct a sutured Heegaard diagram
$(\Sigma,\mbox{\boldmath${\alpha}$},
\mbox{\boldmath$\beta_0$}\cup\{\mu\},w,z)$ from $F$, as in the
proof of \fullref{adjunct}. We can assume this diagram already
satisfies the conclusion of \fullref{support}. The
subsurface $\mathcal P\subset\Sigma$ is more or less a punctured
$F$. We can extend $\varphi$ by identity to a diffeomorphism of
$\Sigma$. $(Y',K')$ has Heegaard diagram
$(\Sigma,\varphi(\mbox{\boldmath${\alpha}$}),
\mbox{\boldmath$\beta_0$}\cup\{\mu\},w,z)$. The generators of
$\widehat{CFK}(Y,K,-g)$ and $\widehat{CFK}(Y',K',-g)$ are the
same. It is easy to see the boundary maps are also the same, since
the boundary of a holomorphic disk does not pass through $\alpha$
curves inside $\mathcal P$. \end{proof}

\begin{proof}[The second proof] Suppose $\gamma$ is a circle in $F$,
with a framing induced by $F$. Since $\gamma$ can be isotoped off
$F$, we have the surgery exact triangle (see
\cite[Theorem 8.2]{OSz2}):
$$\cdots\hspace{-1pt}\to\widehat{HFK}(Y_{-1}(\gamma),K,-g)\to\widehat{HFK}(Y_{0}(\gamma),K,-g)\to\widehat{HFK}(Y,K,-g)\to\hspace{-1pt}\cdots.$$
$(Y_{0}(\gamma),K)$ has a Seifert surface with genus $<g$, which
is obtained by surgering $F$ along $\gamma$. By the adjunction
inequality, $\widehat{HFK}(Y_{0}(\gamma),K,-g)=0$. Hence our
result holds when $\varphi$ is the positive Dehn twist along
$\gamma$. The result also holds when $\varphi$ is a negative Dehn
twist, since a negative Dehn twist is just the inverse of a
positive one. The general case follows since every
self-diffeomorphism of $F$ is the product of Dehn
twists.\end{proof}

\section{Murasugi sum}

In Gabai's theory of sutured manifold decomposation, the longitude
$\lambda$ of a knot often serves as the suture (see \cite{G3}). So
\fullref{support} says that the boundary map of
$\widehat{CFK}(Y,K,-g)$ ``avoids the suture". This justifies the
name ``sutured Heegaard diagram". Using sutured Heegaard diagrams,
we can give a formula for Murasugi sum. This is our first attempt
to apply our method to sutured manifold decomposition.

\begin{defn}\label{MuraDef}
$F^{(k)}$ is an oriented compact surface in the manifold
$Y^{(k)}$, $k=1,2$. $B^{(k)}\subset Y^{(k)}$ is a 3--ball,
$\int(B^{(k)})\cap F^{(k)}=\emptyset$, $\partial B^{(k)}\cap
F^{(k)}$ is a disk $D^{(k)}$, and $D^{(k)}\cap \partial F^{(k)}$
consists of $n$ disjoint arcs. Remove $\int(B^{(k)})$ from
$Y^{(k)}$, glue the punctured $Y^{(1)}$ and punctured $Y^{(2)}$
together by a homeomorphism of the boundaries, so that $D^{(1)}$
is identified with $D^{(2)}$, and $\partial D^{(1)}\cap\partial
F^{(1)}$ is identified with the closure of $\partial
D^{(2)}-\partial F^{(2)}$. We then get a new manifold
$Y=Y^{(1)}\#Y^{(2)}$ and a surface $F=F^{(1)}\cup F^{(2)}$. Then
$F$ is called the {\it Murasugi sum} of $F^{(1)}$ and $F^{(2)}$,
denoted by $F^{(1)}*F^{(2)}$.

$L^{(k)}=\partial F^{(k)}$ is an oriented link in $Y^{(k)}$. Then
$L=\partial F$ is called the {\it Murasugi sum} of $L^{(1)}$ and
$L^{(2)}$, denoted by $L^{(1)}*L^{(2)}$.
\end{defn}

When $n=1$, this operation is merely connected sum; when $n=2$,
this operation is also known as ``plumbing".

Gabai showed that Murasugi sum is a natural geometric operation in
\cite{G1} and \cite{G2}. We summarize some of his results here:

\begin{thm}[Gabai]
With notation as above, we have:

{\rm (i)}\qua $F$ is a Seifert surface with maximal Euler
characteristic for $L$, if and only if $F^{(1)}$ and $F^{(2)}$ are
Seifert surfaces with maximal Euler characteristic for $L^{(1)}$
and $L^{(2)}$, respectively.

{\rm (ii)}\qua $Y-L$ fibers over $S^1$ with fiber $F$, if and only if
$Y^{(k)}-L^{(k)}$ fibers over $S^1$ with fiber $F$ for $k=1,2$.
\end{thm}

Our result about Murasugi sum is an analogue of Gabai's theorem in
the world of knot Floer homology. We first consider the case of
knots.

\begin{prop}\label{MuraKnot}
Suppose knot $(Y,K)$ is the Murasugi sum of two knots
$(Y^{(1)},K^{(1)})$ and $(Y^{(2)},K^{(2)})$. Genera of
$F,F^{(1)},F^{(2)}$ are $g,g^{(1)},g^{(2)}$, respectively. Then
$$\widehat{CFK}(Y,K,-g)\cong\widehat{CFK}(Y^{(1)},K^{(1)},-g^{(1)})\otimes\widehat{CFK}(Y^{(2)},K^{(2)},-g^{(2)})$$
as ungraded chain complexes. In particular, for any field $\mathbb
F$,
$$\widehat{HFK}(Y,K,-g;\mathbb
F)\cong\widehat{HFK}(Y^{(1)},K^{(1)},-g^{(1)};\mathbb
F)\otimes\widehat{HFK}(Y^{(2)},K^{(2)},-g^{(2)};\mathbb F)$$ as
linear spaces.
\end{prop}
\begin{proof} The proof consists of 3 steps. First of all, starting from the surfaces, we
construct a Heegaard splitting for the pair $(Y,K)$. Secondly, we
explicitly give the $\alpha$ and $\beta$ curves on the Heegaard
surface, hence we get a Heegaard diagram. This diagram is a
sutured Heegaard diagram. Finally, based on the diagram, we prove
our desired formula by using \fullref{support}.

\medskip
{\bf Step 1}\qua {\sl Construct a Heegaard splitting for
$(Y,K)$}

$Y$ is separated into two parts $Y^{(1)}-\int(B^{(1)})$ and
$Y^{(2)}-\int(B^{(2)})$. Let $D=F^{(1)}\cap F^{(2)}$ be a
$2n$--gon. Thicken $F$ to $F\times[0,1]$ in $Y$. Add a 1--handle
$\mathcal H$ to connect $D\times0$ to $D\times1$, so that there is
a simple closed curve $\mu\subset\partial
((F\times[0,1])\cup\mathcal H)$, which bounds a disk in
$Y-((F\times[0,1])\cup\mathcal H)$, and passes through $\mathcal
H$ once.

Add $r^{(k)}$ 1--handles $\mathcal H^{(k)}_{2g+2},\dots,\mathcal
H^{(k)}_{2g+1+r^{(k)}}$ in $Y^{(k)}-B^{(k)}$ to $F^{(k)}\times0$,
as when we construct the Heegaard splitting of
$(Y^{(k)},K^{(k)})$. After handlesliding, we can assume each
1--handle is added to a connected component of $(F-D)\times0$. Let
$U_0$ be the handlebody obtained by adding the $1+r^{(1)}+r^{(2)}$
1--handles to $F\times[0,1]$. $Y-\int(U_0)$ is also a handlebody,
we hence get a Heegaard splitting for $(Y,K)$.

\medskip
{\bf Step 2}\qua {\sl Construct a sutured Heegaard diagram}

$F$ is the Murasugi sum of $F^{(1)}$ and $F^{(2)}$, hence
$\chi(F)=\chi(F^{(1)})+\chi(F^{(2)})-1$. Since $K,K^{(1)},K^{(2)}$
are all knots, we have $g=g^{(1)}+g^{(2)}$.

We have the Mayer--Vietoris sequence:
\begin{eqnarray*}
0&\to&\homo_2(F^{(k)},\partial
F^{(k)})\to\oplus^n\homo_1(I,\partial I)
\\&\to&\homo_1(D,D\cap\partial
F^{(k)})\oplus\homo_1(F^{(k)}-D,(F^{(k)}-D)\cap\partial
F^{(k)})\\
&\to&\homo_1(F^{(k)},\partial F^{(k)})\to0.
\end{eqnarray*}
Here $I$ denotes a segment. It follows that
$$\text{rank}\;\homo_1(F^{(k)}-D,(F^{(k)}-D)\cap\partial F^{(k)})=\text{rank}\;\homo_1(F^{(k)},\partial F^{(k)}),$$
and the map
$$\homo_1(F^{(k)}-D,(F^{(k)}-D)\cap\partial F^{(k)})\to\homo_1(F^{(k)},\partial F^{(k)})$$
is injective. Hence we can choose $2g^{(k)}$ disjoint arcs
$\xi^{(k)}_{2},\dots,\xi^{(k)}_{2g^{(k)}+1}\subset F^{(k)}-D$
representing a basis of $\homo_1(F^{(k)},\partial F^{(k)})$. We
choose the arcs such that they do not separate the two feet of any
1--handle added to $(F-D)\times0$. Hence for each 1--handle
$\mathcal H^{(k)}_j$ added to $(F-D)\times0$, there is a simple
closed curve
$$\tau^{(k)}_j\subset\partial(((F-D)\times[0,1])\cup\mathcal H^{(k)}_j)-(F-D)\times1,$$
such that $\tau^{(k)}_j$ passes through $\mathcal H^{(k)}_j$ once,
and is disjoint from all $\xi^{(k)}_i$'s.

Construct a sutured Heegaard diagram
$$(\Sigma^{(k)},\mbox{\boldmath${\alpha}$}^{(k)},
\mbox{\boldmath$\beta_0$}^{(k)}\cup\{\mu^{(k)}\},w^{(k)},z^{(k)})$$
for $(Y^{(k)},K^{(k)})$ as in the proof of \fullref{adjunct}.
Here $\mbox{\boldmath${\alpha}$}^{(k)}$ consists of
$\alpha^{(k)}_i$'s and $\widetilde{\alpha}^{(k)}_j$'s.
$\alpha^{(k)}_i$ is basically the ``double" of $\xi^{(k)}_i$ when
$2\le i\le2g^{(k)}+1$, and $\widetilde{\alpha}^{(k)}_j$ is the
belt circle of $\mathcal H^{(k)}_j$.

We can assume $\alpha^{(k)}_1\subset D$ is the boundary of a
smaller $2n$--gon concentric to $D$, and the restrictions of
$\mbox{\boldmath$\beta_0$}$ curves in $D\times\{0,1\}$ are
segments perpendicular to the corresponding edges of $D$.
Moreover, let the relative positions of $\mu^{(1)}\cap D$ and
$\mu^{(2)}\cap D$ in $D$ be the same. See \fullref{Fig:3} for the case
when $D$ is a square.

Suppose the sides of $D$ are $a_1,b_1,a_2,b_2,\dots,a_n,b_n$ in
cyclic order, where $a_i\subset\partial F^{(1)}$,
$b_i\subset\partial F^{(2)}$. In $\Sigma^{(k)}$, there is a
subsurface $Q^{(k)}$, which is the union of a punctured
$D\times0$, a punctured $D\times1$ and a tube whose belt circle is
$\alpha^{(k)}_1$. We glue $\Sigma^{(1)}-\cup_{i=1}^n
\int(a_i\times[0,1])$ and $\Sigma^{(2)}-\cup_{i=1}^n
\int(b_i\times[0,1])$ together, such that $Q^{(1)}$ is identified
with $Q^{(2)}$, and the edge $(a_i\cap b_{i\pm1})\times[0,1]$ in
$\Sigma^{(1)}$ is glued to the $(a_i\cap b_{i\pm1})\times[0,1]$ in
$\Sigma^{(2)}$.

After the gluing, we get a surface $\Sigma$. We also identify
$\alpha^{(1)}_1$ with $\alpha^{(2)}_1$, $\mu^{(1)}$ with
$\mu^{(2)}$, $(w^{(1)},z^{(1)})$ with $(w^{(2)},z^{(2)})$. The
objects after identification are called $Q$, $\alpha_1$, $\mu$,
$(w,z)$, respectively.

\begin{figure}[ht!]\small\label{Fig:3}
\labellist\hair1pt
\pinlabel $\alpha_1^{(1)}$ [t] <3pt,0pt> at 154 711
\pinlabel $\alpha_1^{(2)}$ [t] <3pt,0pt> at 428 711
\pinlabel {$\mu\cap D$} at 246 613
\endlabellist
\cl{\includegraphics[width=4in]{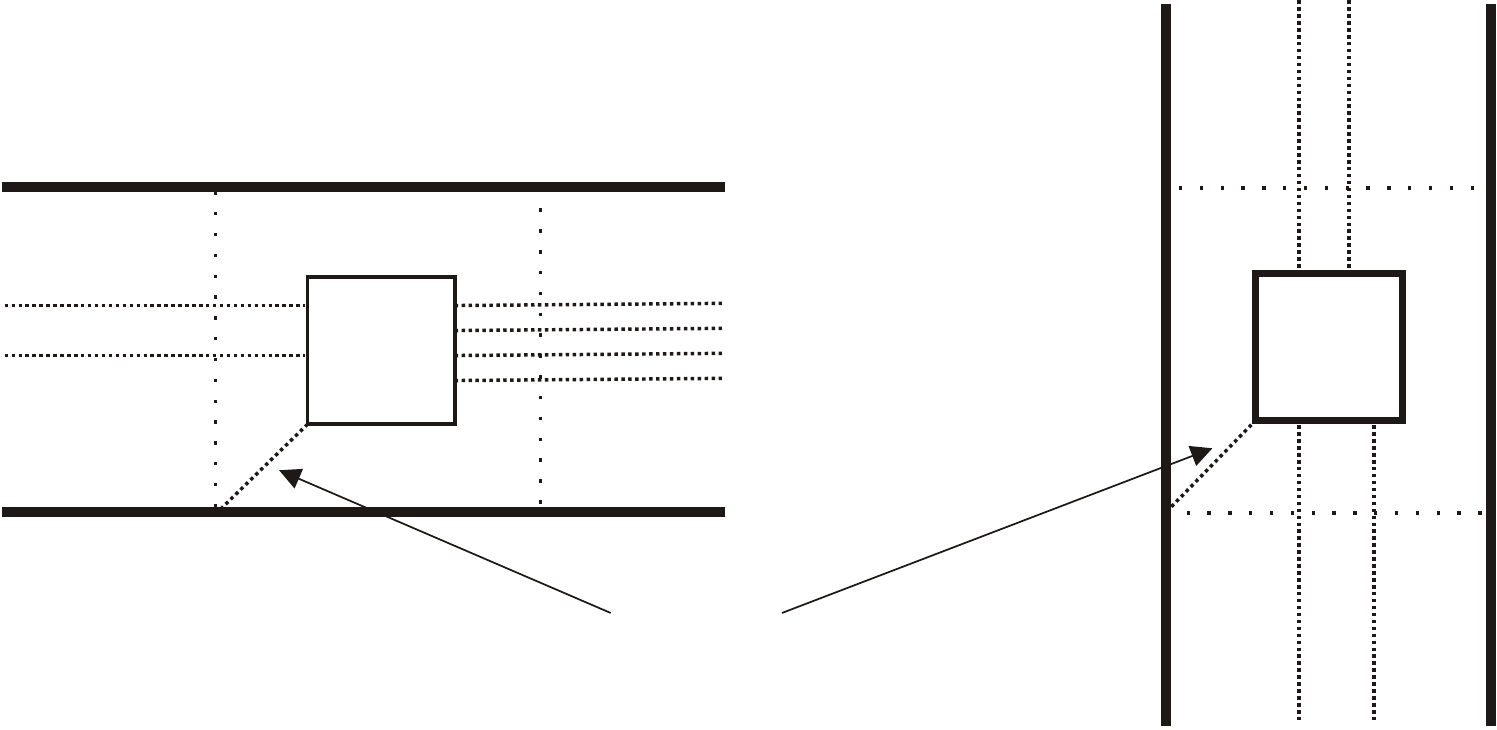}}
\caption{Local pictures of $D^{(k)}\times\{0,1\}$}
\end{figure}

Now we have a diagram
$$(\Sigma,\mbox{\boldmath${\alpha}$}^{(1)}\cup\mbox{\boldmath${\alpha}$}^{(2)},
\mbox{\boldmath$\beta_0$}^{(1)}\sqcup\mbox{\boldmath$\beta_0$}^{(2)}\cup\{\mu\},w,z)$$
We can fit the curves
$\mbox{\boldmath${\alpha}$}^{(1)}\cup\mbox{\boldmath${\alpha}$}^{(2)},
\mbox{\boldmath$\beta_0$}^{(1)}\sqcup\mbox{\boldmath$\beta_0$}^{(2)}\cup\{\mu\}$
into the Heegaard splitting we got in Step 1, so that each curve
bounds a disk in some handlebody. Now it is easy to see
$$(\Sigma,\mbox{\boldmath${\alpha}$}^{(1)}\cup\mbox{\boldmath${\alpha}$}^{(2)},
\mbox{\boldmath$\beta_0$}^{(1)}\sqcup\mbox{\boldmath$\beta_0$}^{(2)}\cup\{\mu\},w,z)$$
is a sutured Heegaard diagram for $(Y,K)$. It is understood that
the ``suture" is the longitude $\lambda$ of $K$.

\medskip
{\bf Step 3}\qua {\sl Prove the formula}

By \fullref{support} and \fullref{support2}, in order
to compute $\widehat{HFK}(Y,K,-g)$, we only need to wind
$\widetilde{\alpha}^{(k)}_j$ along $\tau^{(k)}_j$ many times, and
count the holomorphic disks which are disjoint from $\lambda$.
Suppose $\Phi$ is such a disk, then the local multiplicities of
$\Phi$ at the vertices of $D\times\{0,1\}$ are all 0. As in
\fullref{Fig:4}, we now find that $\Phi\cap Q$ is separated into two
disjoint parts, one is extended into $\Sigma^{(1)}-Q$, the other
is extended into $\Sigma^{(2)}-Q$. (Each part itself may be
disconnected or empty.)

Since $Q$ is the only common part of $\Sigma^{(1)}$ and
$\Sigma^{(2)}$, we now conclude that $\Phi$ consists of two
disjoint parts, one is a holomorphic disk $\Phi^{(1)}$ in
$\Sigma^{(1)}$, the other is a holomorphic disk $\Phi^{(2)}$ in
$\Sigma^{(2)}$. $\Phi^{(k)}$ is a holomorphic disk for
$\widehat{CFK}(Y^{(k)},K^{(k)},-g^{(k)})$.

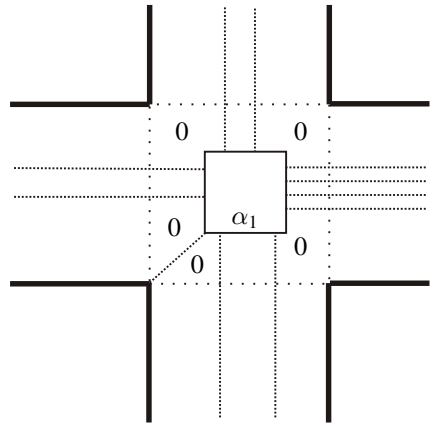
\begin{figure}[ht!]\small\label{Fig:4}
\labellist\hair 2pt
\pinlabel $\alpha_1$ [b] at 292 641
\pinlabel 0 at 321 634
\pinlabel 0 at 321 694
\pinlabel 0 at 259 694
\pinlabel 0 at 255 644
\pinlabel 0 at 267 625
\endlabellist
\cl{\includegraphics[width=2.2in]{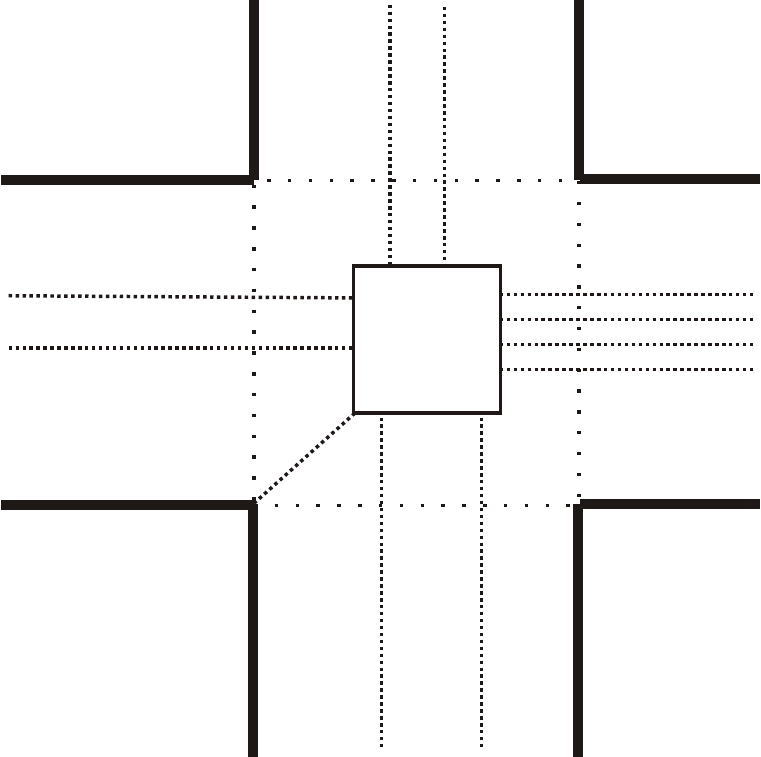}}
\caption{Local pictures of
$D\times\{0,1\}$ after the operation}
\end{figure}

Conversely, if we have holomorphic disks $\Phi^{(k)}$ for
$\widehat{CFK}(Y^{(k)},K^{(k)},-g^{(k)})$, $k=1,2$, then they are
disjoint in $\Sigma$, since they are disjoint from
$\lambda^{(k)}$. Now we can put them together to get a holomorphic
disk $\Phi$ for $\widehat{CFK}(Y,K,-g)$.

Now the formula is obvious.
\end{proof}

Before dealing with the case of links, we recall the definition of
knot Floer homology for links. In \cite{OSz2}, Ozsv\'ath and
Szab\'o gave a well-defined correspondence from links to knots,
and the homology for links is defined to be the homology for the
corresponding knots.

The construction is described as follows: given any
null-homologous oriented $n$--component link $L$ in $Y$, choose
two points $p,q$ on different components of $L$. Remove two balls
at $p,q$, then glue in $S^2\times I$. Inside $S^2\times I$, there
is a band, along which we can perform a connected sum of the two
components of L containing $p$ and $q$. We choose the band so that
the connected sum respects the original orientation on $L$. Now we
have a link in $Y^3\#S^2\times S^1$, with one fewer component.
Repeat this construction until we get a knot. The new knot is
denoted by $\kappa(L)$, and the new manifold is denoted by
$\kappa(Y)=Y\#(|L|-1)(S^2\times S^1)$. Ozsv\'ath and Szab\'o
proved that this correspondence
$(Y,L)\mapsto(\kappa(Y),\kappa(L))$ is well-defined.

Define $\Pi$ to be a link in $S^2\times S^1$, such that $\Pi$
consists of two copies of $\textsl{point}\times S^1$, but with
different orientations. $\Pi$ is a fibred link, its fiber is an
annulus. When we do plumbing of $\Pi$ with other links, we always
choose the annulus as the Seifert surface for $\Pi$. It is not
hard to see that Ozsv\'ath and Szab\'o's construction is more or
less doing plumbing with copies of $\Pi$. (See \cite{Ni} for an
explanation.)

The next lemma is a special case of our general theorem about
Murasugi sum.

\begin{lem}\label{plumbing}
Suppose $(Y,L)$ is a link with Seifert surface $F$. Do plumbing
for $L$ and $\Pi$, we get a link $(Y'=Y\#S^2\times S^1,L')$ with
Seifert surface $F'$. Then
$$\widehat{HFK}(Y,L,\frac{|L|-\chi(F)}2)\cong\widehat{HFK}(Y',L',\frac{|L'|-\chi(F')}2)$$
as abelian groups.
\end{lem}
\begin{proof}
If the plumbing merges two components of $L$, then the result
holds by the discussion before this lemma. Now we consider the
case that the plumbing splits a component of $L$ into two
components. Without loss of generality, we can assume $L$ is a
knot.

Now $L'=L*\Pi$ is a two-component link, we have to consider the
knot $\kappa(L')$. But $\kappa(L')$ is just the plumbing of $L$
with the fibred knot $\Pi*\Pi$. Hence the result holds by
\fullref{MuraKnot}.
\end{proof}

Now we can give our general theorem about Murasugi sum.

\begin{thm}
Notations as in \fullref{MuraDef}. Let $\mathfrak
i(F)=\frac{|\partial F|-\chi(F)}2$. Given a field $\mathbb F$, we
have
$$
\widehat{HFK}(Y,L,\mathfrak i(F))\cong
\widehat{HFK}(Y^{(1)},L^{(1)},\mathfrak
i(F^{(1)}))\otimes\widehat{HFK}(Y^{(2)},L^{(2)},\mathfrak
i(F^{(2)}))
$$
as linear spaces. Here we use $\mathbb F$--coefficients.
\end{thm}
\begin{proof}
Suppose the Murasugi sum is done along a $2n$--gon $D$. The sides
of $D$ are denoted by $a_1,b_1,a_2,b_2,\dots,a_n,b_n$ in cyclic
order, where $a_i\subset\partial F^{(1)}$, $b_i\subset\partial
F^{(2)}$.

Push a neighborhood of $b_i$ slightly out of $D$, to get a
rectangle $R(b_i)$. We do plumbing of $F^{(1)}$ with $n-1$ copies
of $\Pi$, along $R(b_1),R(b_2),\dots,R(b_{n-1})$. We get a new
link $L^{(1)}_1$ with Seifert surface $F^{(1)}_1$. There is an arc
$a\subset L^{(1)}_1$, $a$ contains $a_1,\dots,a_n$ in order. \fullref{plumbing} shows that
$$\widehat{HFK}(L^{(1)}_1,\mathfrak i(F^{(1)}_1))\cong\widehat{HFK}(L^{(1)},\mathfrak i(F^{(1)})).$$
(We suppress the ambient 3--manifolds in the formula.)

Similarly, we plumb $F^{(2)}$ with $n-1$ copies of $\Pi$, along
$R(a_2),R(a_3),\dots,R(a_n)$, to get a link $L^{(2)}_1$. There is
an arc $b\subset L^{(2)}_1$, $b$ contains $b_1,\dots,b_n$ in
order. Moreover,
$$\widehat{HFK}(L^{(2)}_1,\mathfrak i(F^{(2)}_1))\cong\widehat{HFK}(L^{(2)},\mathfrak i(F^{(2)})).$$
We perform the operation $\kappa$ to $L^{(1)}_1, L^{(2)}_1$, with
all the connecting bands added outside $a,b$, to get new knots
$K^{(1)},K^{(2)}$. $K^{(1)}$ contains $a_1,\dots,a_n$ in cyclic
order, and $K^{(2)}$ contains $b_1,\dots,b_n$ in cyclic order.

Now it is easy to see the Murasugi sum $K=K^{(1)}*K^{(2)}$ is
still a knot. Our result holds by \fullref{MuraKnot} and
\fullref{plumbing}.
\end{proof}

\section{Semifibred satellite knots}

The reader should note that some notations in this section are
different from the last section. This is not very disturbing,
since this section is independent of the last one.

\begin{defn}\label{satellite}
Suppose $K$ is a null-homologous knot in $Y$, $F$ is a Seifert
surface of $K$ (not necessarily has minimal genus). $V$ is a
3--manifold, $\partial V=T^2$, $L\subset V$ is a nontrivial knot.
$G\subset V$ is a compact connected oriented surface so that $L$
is a component of $\partial G$, and $\partial G-L$ (may be empty)
consists of parallel essential circles on $\partial V$.
Orientations on these circles are induced from the orientation on
$G$, we require that these circles are parallel as oriented ones.
We glue $V$ to $Y-\int(N(K))$, so that any component of $\partial
G-L$ is null-homologous in $Y-\int(N(K))$. The new manifold is
denoted by $Y^*$, and the image of $L$ in $Y^*$ is denoted by
$K^*$. We then say $K^*$ is a {\it satellite knot} of $K$, and $K$
a {\it companion knot} of $K^*$. Let $p$ denote the number of
components of $\partial G-L$, $p$ will be called the {\it winding
number} of $L$ in $V$.

Moveover, if $V-L$ fibers over the circle so that $G$ is a fiber
and $\chi(G)<0$, then we say $K^*$ is a {\it semifibred satellite
knot}.
\end{defn}

\begin{rem}
In order to avoid some trivial cases, we often need some
additional condition on $G$ in the definition of satellite knot,
say, incompressible in $V-L$. But this would not affect the
results stated in this paper.
\end{rem}

The classical case is $Y=Y^*=S^3$, and $V$ is a solid torus. A
large number of classical satellite knots are semifibred. For
example, it is well-known that cable knots are semifibred, see
\cite[10.I]{Ro}. A bit more work can show that if $L$ is a
``homogeneous braid" in the solid torus $V$, then $K^*$ is
semifibred (see \cite{St}).

Our goal in this section is
\begin{thm}\label{Gsemifibred}
Notations as in \fullref{satellite}. $K^*$ is a semifibred
satellite knot. Suppose the genera of $F,G$ are $g,h$,
respectively, and the winding number is $p$. Then
$$\widehat{HFK}(Y^*,K^*,-(pg+h))\cong\widehat{HFK}(Y,K,-g)$$
as abelian groups.
\end{thm}

In the case of classical knots, our result should be compared with
the well-known relation for Alexander polynomial:
$$\Delta_{K^*}(t)=\Delta_K(t^p)\Delta_{L}(t).$$
See, for example, \cite{BZ}.

\begin{rem}
Not many results were known previously on the knot Floer homology
of satellite knots. Eftekhary computed the top filtration term for
Whitehead doubles in \cite{Ef}. And some terms for $(p,pn\pm1)$
cable knots were computed by Hedden in \cite{He}. Hedden and
Ording also have an ongoing program to compute the Floer homology
of $(1,1)$ satellite knots. Whitehead doubles are not semifibred
in our sense, although $V-L$ fibers over $S^1$ in this case.
\end{rem}

\begin{construction}
As the reader may have found in the last section, we have to spend
most efforts on the description of the construction of a suitable
Heegaard diagram, although the idea of such construction is very
simple. Our construction here consists of 5 steps. The notations
are as before. In this construction, we assume the monodromy of
the fibred part $V-\int(N(L))$ is a special map $\psi$, which will
be defined in Step 2.

\medskip{\bf Step 0}\qua {\sl A Heegaard splitting of $Y^*$}

A Heegaard splitting of $(Y^*,K^*)$ can be constructed as follows.
Pick $p$ parallel copies of $F$: $F^{(1)},\dots,F^{(p)}$. Glue
them to $G$, so as to get a surface $F^*$ of genus $pg+h$. Thicken
$F^*$ to $F^*\times[0,1]$ in $Y^*$. Add a one handle $\mathcal
H^{(k)}$ connecting $F^{(k)}\times1$ and $F^{(k+1)}\times0$,
$k=1,\dots,p$, here $F^{(p+1)}=F^{(1)}$. Add a 1--handle $\mathcal
H^*$ connecting $G\times1$ to $G\times0$, so that it is parallel
to $\mathcal H^{(k)}$. Then add $r$ 1--handles to $F^{(1)}\times0$
in the same way as when we constructed the Heegaard splitting of
$(Y,K)$ in the proof of \mbox{\fullref{adjunct}}. See \fullref{Fig:7} for a
schematic picture.

Now we have a Heegaard splitting $Y^*=U_0^*\cup U_1^*$, $U_0^*$ is
the union of $F^*\times[0,1]$ and some 1--handles. In the rest of
this construction, we will construct the corresponding Heegaard
surface abstractly, and give the $\alpha$ and $\beta$ curves on
this Heegaard surface. Hence we get a Heegaard diagram which can
be fit into the Heegaard splitting we construct in Step 0.

\medskip{\bf Step 1}\qua {\sl Construct block surfaces with curves
on them}

$A$ is a genus $g$ surface with boundary consisting of two
circles, denoted by $\alpha_1,\lambda$. $A$ is basically a
punctured $F\times1$. Pick an arc $\delta$ connecting $\alpha_1$
to $\lambda$. Then we can choose two $2g$--tuples of mutually
disjoint proper arcs in $A-\delta$: $(\xi_2,\dots,\xi_{2g+1})$ and
$(\eta_2,\dots,\eta_{2g+1})$, so that
$\partial\xi_i\subset\lambda$, $\partial\eta_i\subset\alpha_1$.
Moreover, $\eta_i$ is disjoint from $\xi_j$ when $j\ne i$, and
$\eta_i$ intersects $\xi_i$ transversely in a single intersection
point.  The reader is referred to \cite[Figure~1]{OSz3} for a
choice of these curves.

Let $B^*$ be a genus $h$ surface with boundary consisting of two
circles, denoted by $\alpha_1^*,\lambda^*$. We can choose curves
$\xi^*_j,\eta^*_j,\delta^*$ as before. For the arc $\delta^*$
connecting $\alpha^*_1$ to $\lambda^*$, we pick $p$ parallel
copies, $\delta_1,\dots,\delta_p$, lying on the same side of
$\delta^*$. Choose a point in each $\delta_k$, remove a small open
disk at each chosen point, then get a surface with boundary
consisting of $p+2$ circles
$\alpha_1^*,\lambda^*,\lambda^{(1)},\dots,\lambda^{(p)}$, called
$B$. The remaining part of $\delta_k$ consists of two arcs
$\pi^{(k)},\rho^{(k)}$, here $\pi^{(k)}$ connects $\lambda^*$ to
$\lambda^{(k)}$. See \mbox{\fullref{Fig:5}} for the local picture.

\begin{figure}[ht!]\small\label{Fig:5}
\labellist\hair 2pt
\pinlabel $\pi^{(3)}$ [t] at 228 645
\pinlabel $\pi^{(2)}$ [t] at 228 610
\pinlabel $\pi^{(1)}$ [t] at 228 565 
\pinlabel $\lambda^{(3)}$ <1pt,2pt> at 303 653
\pinlabel $\lambda^{(2)}$ <1pt,2pt> at 303 610
\pinlabel $\lambda^{(1)}$ <1pt,2pt> at 303 566
\pinlabel $\delta^*$ [t] at 303 523
\pinlabel $\rho^{(1)}$ [t] at 377 566
\pinlabel $\rho^{(2)}$ [t] at 377 610
\pinlabel $\rho^{(3)}$ [t] at 377 645
\pinlabel $\lambda`^*$ [b] <0pt,1pt> at 139 507
\pinlabel $\alpha_1^*$ [b] at 466 506
\pinlabel $\xi^*_j$ [b] at 30 609
\pinlabel* $\eta^*_j$ [tl] at 553 627
\endlabellist
\cl{\includegraphics[width=4.5in]{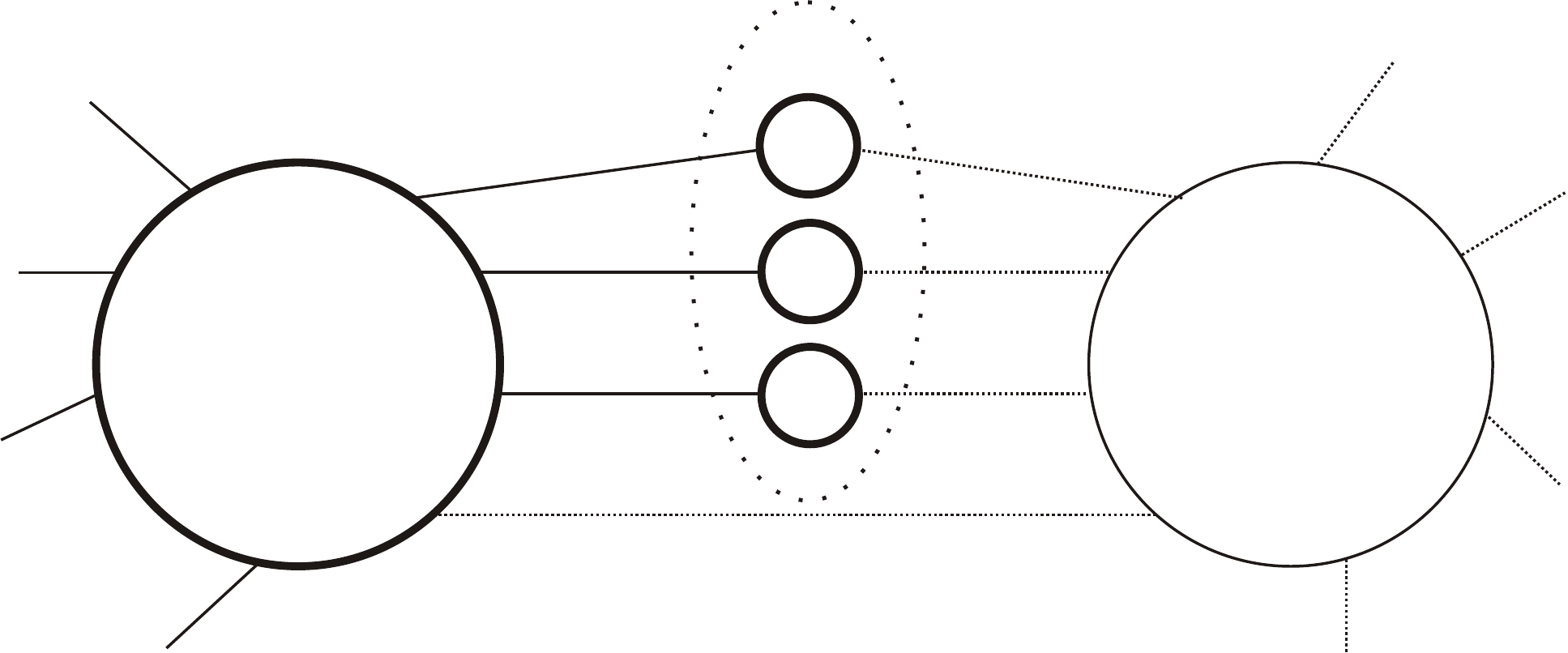}}
\caption{Local picture near $\delta^*$}
\end{figure}
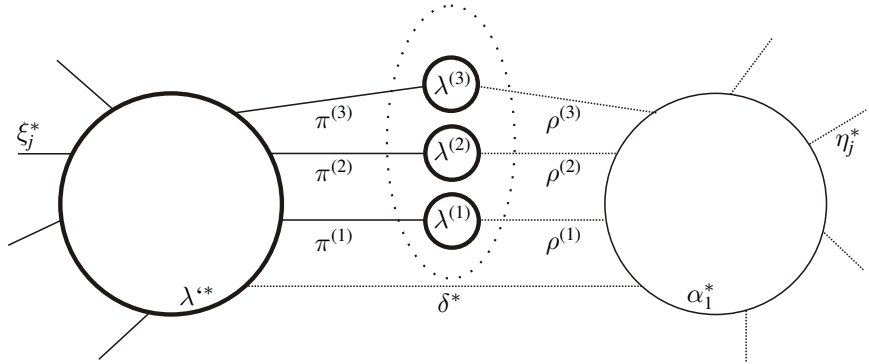

\medskip{\bf Step 2}\qua {\sl Construct the monodromy $\psi$}

Take $G\times[0,2]$, glue the two ends together by the identity,
so as to get $G\times S^1$. Two surfaces $G\times0,G\times1$
separate $G\times S^1$ into two parts. Choose a small disk $D_1$
in the interior of $G$. Denote $G\times[0,1]\cup D_1\times[1,2]$
by $V_0$, $G\times[1,2]-\int(D_1)\times[1,2]$ by $V_1$.

$\overline B$ denotes the copy of $B$ reflected across its
boundary. Curves on $\overline B$ are denoted by
$\overline{\xi}_j^*$, etc. Glue $B$ and $\overline B$ so that
$\lambda^*,\alpha_1^*$ are identified with
$\overline{\lambda}^*,\overline{\alpha}_1^*$. $B\cup\overline B$
can be naturally identified with the surface
$$((G-\int(D_1))\times1)\cup ((G-\int(D_1))\times0)\cup(\partial
D_1\times[1,2])\cup(\lambda^*\times[0,1]),$$ where $\lambda^*$ is
the boundary component of $G$ which corresponds to the longitude
of $K^*$, by abuse of notation. $B\supset (G-\int(D_1))\times1$,
$\overline B\supset(G-\int(D_1))\times0$.

Glue $\eta_j^*$ and $\overline{\eta}_j^*$ together to a closed
curve $\beta_j^*$, glue $\xi_j^*,\overline{\xi}_j^*$ together to a
closed curve $\alpha_j^*$, $j=2,\dots,2h+1$. Glue $\delta^*$ and
$\overline{\delta}^*$ together to a closed curve $\mu^*$. Glue
$\pi^{(k)}$ and $\overline{\pi}^{(k)}$ together to an
arc\footnotemark. Glue $\rho^{(k)}$ and $\overline{\rho}^{(k)}$
together to an arc. We have the properties:

\footnotetext{$\pi^{(k)}$ connects $\lambda^*$ to $\lambda^{(k)}$,
$\overline{\pi}^{(k)}$ connects $\overline{\lambda}^*$ to
$\overline{\lambda}^{(k)}$. We glued $\lambda^*\subset B$ to
$\overline{\lambda}^*\subset \overline B$, but did not glue
$\lambda^{(k)}$ to $\overline{\lambda}^{(k)}$, so
$\pi^{(k)}\cup\overline{\pi}^{(k)}$ is an arc.}

($\alpha$)\qua The circles $\alpha_j^*$ ($j=1,2,\dots,2h+1$) bound
disks in $V_0$. The arc $\pi^{(k)}\cup\overline{\pi}^{(k)}$
cobounds a half-disk\footnotemark\; in $V_0$ with a vertical arc
on $\lambda^{(k)}\times[0,1]$, $k=1,\dots,p$.

\footnotetext{Of course a half-disk is homeomorphic to a disk. We
use the term ``half-disk" because this disk will be part of a disk
bounded by an $\alpha$ curve constructed later.}

($\beta$)\qua The circles $\mu^*,\beta_j^*$ ($j=2,\dots,2h+1$) bound
disks in $V_1$. The arc $\rho^{(k)}\cup\overline{\rho}^{(k)}$
cobounds a half-disk in $V_1$ with a vertical arc on
$\lambda^{(k)}\times[1,2]$, $k=1,\dots,p$.

As in \fullref{Fig:5}, the dotted circle encloses a $p$--punctured disk
$D_2$. There is a diffeomorphism $\psi$ of $B$, supported in
$D_2$, sending $\lambda^{(k)}$ to $\lambda^{(k-1)}$
($\lambda^{(0)}=\lambda^{(p)}$). We draw $\psi(\rho^{(k)})$ in
\fullref{Fig:6}.

$\psi$ can be extended by identity to a diffeomorphism of $G$,
still denoted by $\psi$. Cut $G\times S^1$ open along $G\times1$,
reglue by $\psi$, namely, glue each point $x\in G\times(1+0)$ to
$\psi(x)\in G\times(1-0)$. Now we get a surface bundle over the
circle, with monodromy $\psi$.  This surface bundle will serve as
our $V-\int(N(L))$.

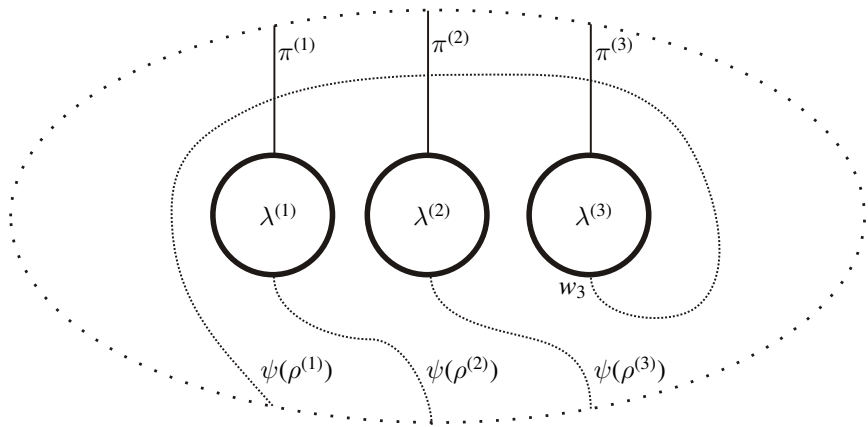
\begin{figure}[ht!]\small\label{Fig:6}
\labellist\hair 2pt
\pinlabel $\pi^{(3)}$ [l] at 371 583
\pinlabel $\pi^{(2)}$ [l] at 287 587
\pinlabel $\pi^{(1)}$ [l] at 208 583
\pinlabel $\lambda^{(3)}$ <3pt,0pt> at 369 498
\pinlabel $\lambda^{(2)}$ <3pt,0pt> at 286 498
\pinlabel $\lambda^{(1)}$ <3pt,0pt> at 207 498
\pinlabel $\psi(\rho^{(1)})$ [bl] at 200 409
\pinlabel $\psi(\rho^{(2)})$ [bl] at 285 409
\pinlabel $\psi(\rho^{(3)})$ [bl] at 371 409
\pinlabel $w_3$ [tr] <1pt,0pt> at 370 465
\endlabellist
\cl{\includegraphics[width=4.5in]{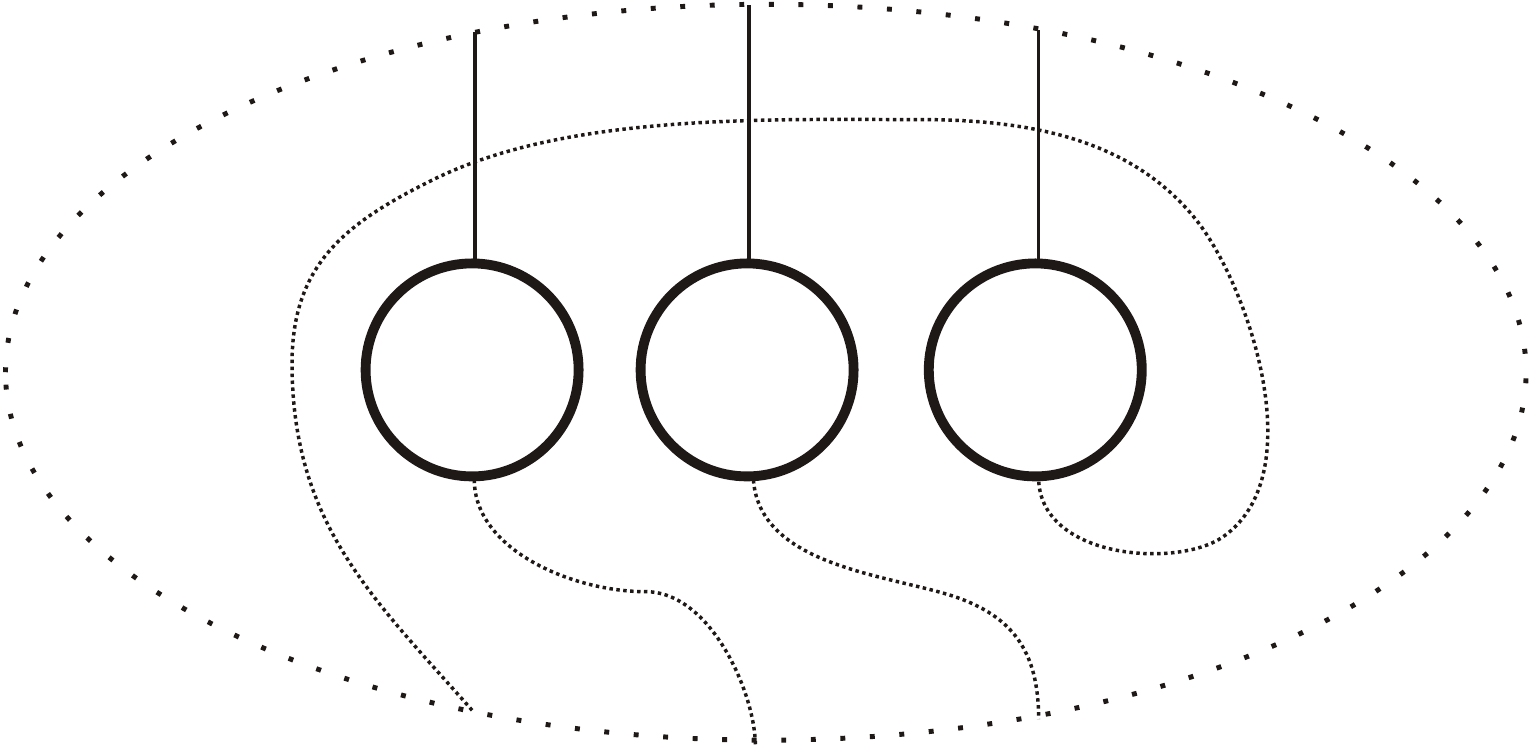}}
\caption{Local picture inside the circle}
\end{figure}

$B\cup \overline B$ naturally lies in $V-\int(N(L))$ as before, and
it separates the bundle into $V_0$ and $V_1$. We will view
$B\cup\overline B$ as living on the boundary of $V_0$. We again
have the curves on $B\cup \overline B$ satisfying Properties
($\alpha$),($\beta$), except that in Property ($\beta$),
$\psi(\rho^{(k)})\cup\overline{\rho}^{(k)}$ cobounds a half-disk
in $V_1$ with a vertical arc on $\lambda^{(k)}\times[1,2]$,
$k=1,\dots,p$. $\partial V=T^2$ is the union of the annuli
$$\lambda^{(1)}\times[0,1],\lambda^{(2)}\times[1,2],\lambda^{(2)}\times[0,1],\dots,\lambda^{(p)}\times[0,1],\lambda^{(1)}\times[1,2].$$

\medskip{\bf Step 3}\qua {\sl Glue blocks together}

Take $p$ copies of $A$: $A^{(1)},A^{(2)},\dots,A^{(p)}$, and let
$\overline A^{(k)}$ denote the copy of $A^{(k)}$ reflected across
its boundary. Curves on $A^{(1)}$ are denoted by
$\xi_i^{(1)},\eta_i^{(1)}$, etc. One may worry about the $\lambda$
curve on $A^{(k)}$, which will be called $\lambda^{(k)}$ by our
convention, and this name coincides with the boundary curve
$\lambda^{(k)}$ of $B$. But since we are going to identify these
two curves, we do not introduce any new notation to distinguish
them. Similarly, we name the curves on $\overline A^{(k)}$ by
$\overline{\xi}_i^{(k)}$, etc.

\begin{figure}[ht!]\small\label{Fig:7}
\labellist
\pinlabel $A^{(3)}$ [b] <5pt,0pt> at 54 761
\pinlabel $A^{(1)}$ [b] <5pt,0pt> at 249 761
\pinlabel $A^{(2)}$ [b] <5pt,0pt> at 442 761
\pinlabel $\wtilde{A}$ [b] <1pt,0pt> at 146 761
\pinlabel $\wbar{A}^{(3)}$ [b] <5pt,0pt> at 537 761
\pinlabel $\wbar{A}^{(2)}$ [b] <5pt,0pt> at 352 761
\pinlabel $\wbar{B}$ [l] at 518 384
\pinlabel $B$ [r] at 43 384
\pinlabel $\psi(\rho^{(3)})$ [t] <3pt,2pt> at 220 451 
\pinlabel $\wbar{\rho}^{(3)}$ [tl] <0pt,4pt> at 497 462
\endlabellist
\vspace{2mm}\cl{\includegraphics[width=4.5in]{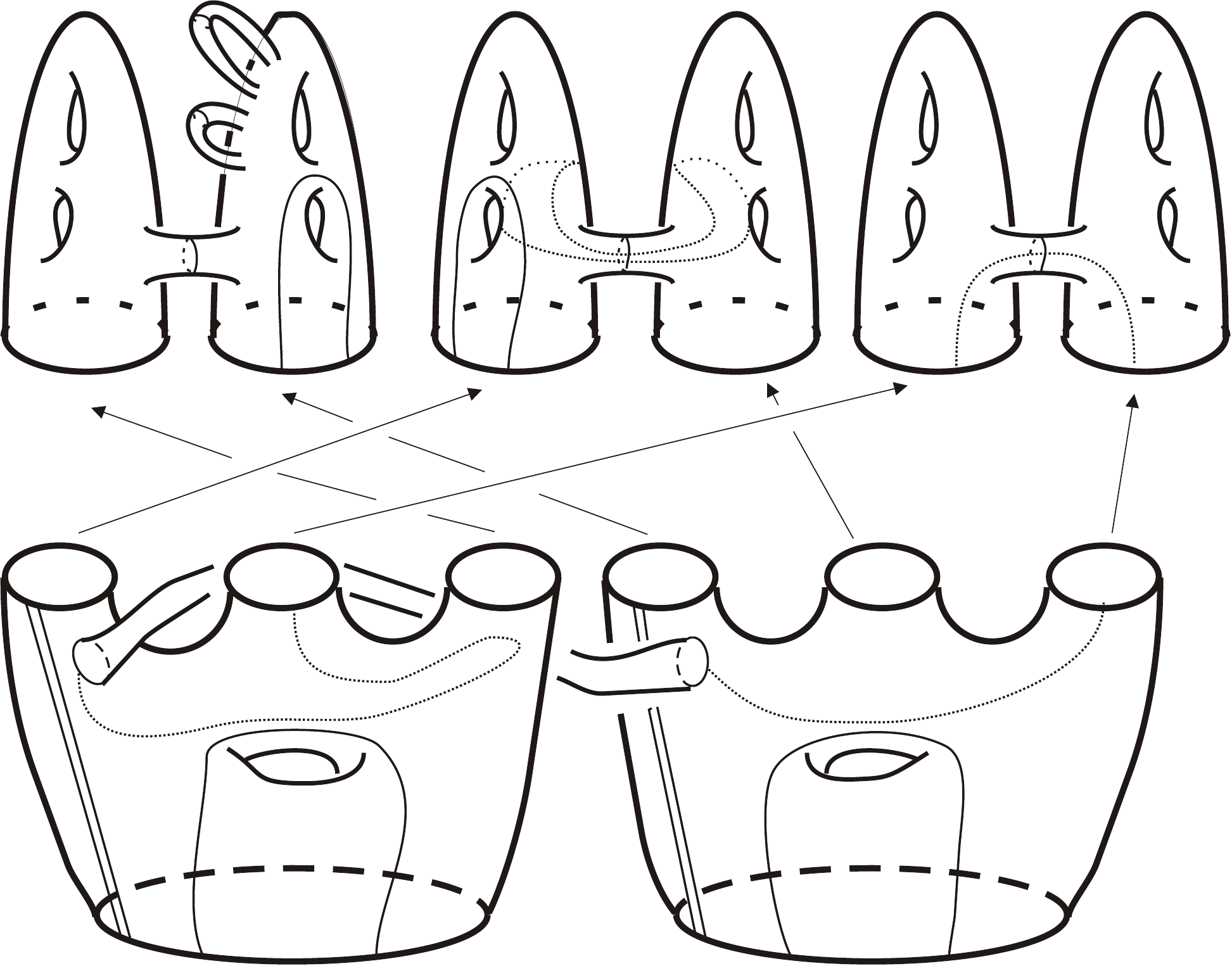}}
\caption{The Heegaard surface $\Sigma^*$
with some sample curves drawn}
\end{figure}
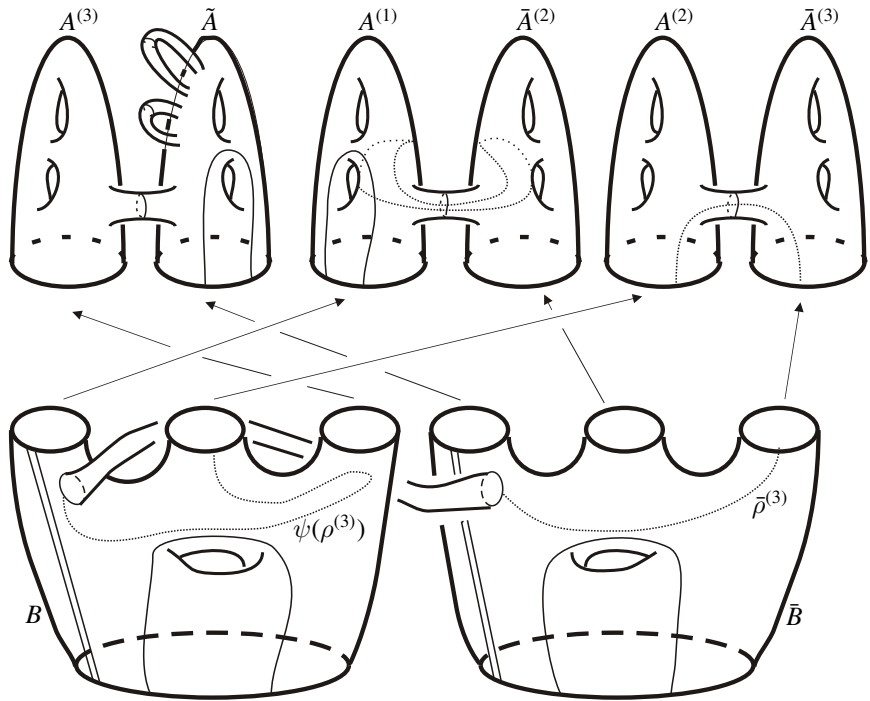

We glue $A^{(k)}$ and $\overline A^{(k+1)}$ so that
$\alpha_1^{(k)}$ is identified with $\overline{\alpha}_1^{(k+1)}$,
glue $A^{(k)}$ and $B$ along $\lambda^{(k)}$, glue $\overline
A^{(k)}$ and $\overline B$ along $\overline{\lambda}^{(k)}$.

Add $r$ tubes to $\overline A^{(1)}$, as in Step 0.  The union of the
$(2r)$--punctured $\overline A^{(1)}$ and the $r$ new tubes is
called $\widetilde{A}$. Meridians of these new tubes are called
$\widetilde{\alpha}_{2g+2},\dots,$ $\widetilde{\alpha}_{2g+1+r}$.
One can choose circles
$\widetilde{\beta}_2,\dots,\widetilde{\beta}_{2g+1+r}$ on $
A^{(p)}\cup\widetilde A$ as in the sutured Heegaard diagram
$(\Sigma,\mbox{\boldmath${\alpha}$},
\mbox{\boldmath$\beta_0$}\cup\{\mu\},w,z)$, so that they are
disjoint from $\delta^{(p)},\overline{\delta}^{(1)}$.

Glue $\eta_i^{(k)}$ and $\overline{\eta}_i^{(k+1)}$ together to a
closed curve $\beta_i^{(k)}$, $k=1,\dots,p-1$,\break $i=2,\dots,2g+1$.
Glue $\xi_i^{(k)},\overline{\xi}_i^{(k)}$,two copies of
$\pi^{(k)}$ and two copies of $\overline{\pi}^{(k)}$ together to a
closed curve $\alpha_i^{(k)}$, $k=1,2,\dots,p,\;i=2,\dots,2g+1$.

Glue $\delta^{(k)},\overline{\delta}^{(k+1)}$,
$\psi(\rho^{(k+1)})$ and
 $\overline{\rho}^{(k+1)}$ together to a closed curve
$\omega^{(k)}$, $k=1,\dots,p$.

\medskip{\bf Step 4}\qua {\sl A Heegaard diagram for $(Y^*,K^*)$}

From last step, we have a surface
$$\Sigma^*=\widetilde A\cup A^{(1)}\cup\overline A^{(2)}\cup A^{(2)}\cup\dots\cup\overline A^{(p)}\cup A^{(p)}\bigcup \overline B\cup B,$$
with two collections of disjoint closed curves
\begin{eqnarray*}
\mbox{\boldmath${\alpha}$}^* &= &
\{\alpha_i^{(k)},k=1,\dots,p,\;i=2,\dots,2g+1\}\bigcup\{\widetilde{\alpha}_{2g+2},\dots,\widetilde{\alpha}_{2g+1+r}\}\\
& &\bigcup
\{\alpha_2^*,\dots,\alpha_{2h+1}^*\}\bigcup\{\alpha_1^*,\alpha_1^{(1)},\dots,\alpha_1^{(p)}\},\\
\mbox{\boldmath${\beta}$}^* &= &
\{\beta_i^{(k)},k=1,\dots,p-1,\;i=2,\dots,2g+1\}\bigcup\{\widetilde{\beta}_{2},\dots,\widetilde{\beta}_{2g+1+r}\}\\
& &\bigcup
\{\beta_2^*,\dots,\beta_{2h+1}^*\}\bigcup\{\omega^{(1)},\dots,\omega^{(p)}\}\cup\{\mu^*\}.
\end{eqnarray*}
We claim that
$(\Sigma^*,\mbox{\boldmath${\alpha}$}^*,\mbox{\boldmath${\beta}$}^*)$
is a Heegaard diagram for $Y^*$. In fact we can fit $\Sigma^*$
into the construction in Step 0, so that $\overline A^{(k)}$ is
basically $F^{(k)}\times0$ with a hole, $A^{(k)}$ is basically
$F^{(k)}\times1$ with a hole, $\overline B$ is basically
$G\times0$ with a hole, $B$ is basically $G\times1$ with a hole.
Now the proof of our claim is as easy as ABC, it is easy to check
the following:

\noindent (A)\qua
$\textrm{genus}(\Sigma^*)=|\mbox{\boldmath${\alpha}$}^*|=|\mbox{\boldmath${\beta}$}^*|=2pg+2h+r+p+1$.

\noindent (B)\qua Curves in $\mbox{\boldmath${\alpha}$}^*$ bound disks
in $U_0^*$, curves in $\mbox{\boldmath${\beta}$}^*$ bound disks in
$U_1^*$.

(For example, in order to check that $\alpha_i^{(k)}$ bounds a
disk in $U_0^*$, we recall that $\alpha_i^{(k)}$ is the union of
$\xi_i^{(k)},\overline{\xi}_i^{(k)}$,two copies of $\pi^{(k)}$ and
two copies of $\overline{\pi}^{(k)}$. $\xi_i^{(k)}$ and
$\overline{\xi}_i^{(k)}$ are two parallel sides of a rectangle
between $A^{(k)}$ and $\overline A^{(k)}$,
$\pi^{(k)}\cup\overline{\pi}^{(k)}$ cobounds a half-disk in $V_0$
with a vertical arc on $\lambda^{(k)}\times[0,1]$, (see Property
($\alpha$) in Step 2,) the union of the rectangle and two copies
of the half-disk is a disk bounded by $\alpha_i^{(k)}$ in
$U_0^*$.)

\noindent (C)\qua $\Sigma^*-\mbox{\boldmath${\alpha}$}^*$ is
connected, $\Sigma^*-\mbox{\boldmath${\beta}$}^*$ is connected.

Pick two points $w^*,z^*$ near $\lambda^*\cap\mu^*$ as understood.
Then
$(\Sigma^*,\mbox{\boldmath${\alpha}$}^*,\mbox{\boldmath${\beta}$}^*,w^*,z^*)$
is a double pointed diagram for $(Y^*,K^*)$.\hfill\qedsymbol
\end{construction}

Strictly speaking, the Heegaard diagram constructed above is not a
sutured Heegaard diagram. In \fullref{suturedHD}, in a
sutured Heegaard diagram, there is a subsurface $\mathcal P$
bounded by 2 curves $\lambda$ and $\alpha_1$. In our diagram here,
the corresponding subsurface
$$\mathcal P^*=A^{(1)}\cup\dots\cup A^{(p)}\cup B$$
has $p+2$ boundary components
$\lambda^*,\alpha^*_1,\alpha^{(1)}_1,\dots,\alpha^{(p)}_1.$

However, we can still handle this diagram by the same method we
used in Section 3. (We could have extended \fullref{suturedHD} to the case of more boundary components, but we
would rather choose the current version for simplicity.)

We claim that the generators of $\widehat{CFK}(Y^*,K^*,-(pg+h))$
lie outside the interior of $\mathcal P^*$. In fact $\chi(\mathcal
P^*)=-2(pg+h)-p$, and we have to choose points on
$\alpha_1^{(1)},\dots,\alpha_1^{(p)}$, where the local
multiplicity is $\frac12$. Argue as in Step 3 of \fullref{adjunct}, we can prove the claim. (We cheat a little bit
here, since the diagram is not known to be weakly admissible now.)

By the construction of the Heegaard diagram for $(Y^*,K^*)$, there
is only one choice for the generators of
$\widehat{CFK}(Y^*,K^*,-(pg+h))$ outside $\widetilde A$. Hence the
domains of the holomorphic disks corresponding to the boundary map
will restrict to relative periodic domains outside of
$\widetilde{A}$. We want to show that these holomorphic disks are
supported inside $A^{(p)}\cup\widetilde A$, hence they are in
one-to-one correspondence with the holomorphic disks for
$\widehat{CFK}(Y,K,-g)$. The basic method is also winding.

\begin{proof}[Proof of \fullref{Gsemifibred}] Suppose the
monodromy of $V-\int(N(L))$ is $\varphi$. Since the fiber $G$ is
connected, and has $p$ parallel components on $\partial V$,
$\varphi$ must permute the $p$ components cyclically. Without loss
of generality, we can assume $\varphi$ sends $\lambda^{(k)}$ to
$\lambda^{(k-1)}$. Hence $\varphi\circ\psi^{-1}$ is isotopic to a
diffeomorphism of $G$, which restricts to identity on $\partial
G$. Here $\psi$ is the monodromy constructed in Step 2 of the
previous construction. Now we change the monodromy of the fibred
part to $\psi$, so as to get a new knot in a new manifold. The new
pair is still denoted by $(Y^*,K^*)$. \fullref{reglue}
says that we only need to prove our theorem for this new knot. The
construction before gives a Heegaard diagram for $(Y^*,K^*)$.

In $B$, we can choose $2h$ simple closed curves
$\tau^*_2,\dots,\tau^*_{2h+1}$, such that they are disjoint from
$\delta^*,\rho^{(k)},\pi^{(k)},\xi^*_j$, except that $\tau^*_j$
intersects $\xi^*_j$ transversely in a single intersection point.

Choose an arc $a\subset\lambda^*$, such that $a$ intersects
$\delta^*$ and all $\xi^*_j$'s, but $a$ is disjoint from
$\pi^{(1)},\dots,\pi^{(p)}$. Let $\mathcal U$ be a small
neighborhood of $a$. Wind $\xi^*_j$'s along $\tau^*_j$'s
sufficiently many times, and apply \fullref{winding}, we find
that if $\mathcal D^*$ is a nonnegative relative periodic domain
in $B$, then $\partial\mathcal D^*$ does not pass through
$\xi^*_j$'s. Moreover, the local multiplicity of $\mathcal D^*$ is
0 at a point $w_p\in\lambda^{(p)}$ near
$\psi(\rho^{(1)})\cap\lambda^{(p)}$ (see \fullref{Fig:6}).

Choose an arc $b\subset\lambda^{(p)}$, $b$ intersects $\pi^{(p)}$,
$b\ni w_p$, but $b$ is disjoint from $\psi(\rho^{(1)})$. Let
$\mathcal V$ be a neighborhood of $b$, $w_p$ is the base point.
Apply \fullref{winding} to the surface $A^{(p)}$, we get the
following conclusion: if $\mathcal D^{(p)}$ is a nonnegative
relative periodic domain in $A^{(p)}$, then $\partial\mathcal
D^{(p)}$ does not pass through $\xi^{(p)}_i$'s. Moreover,
$\partial\mathcal D^{(p)}$ does not pass through $\delta^{(p)}$,
since the $\beta$ curves are disjoint from $\lambda^{(p)}$ while
$\delta^{(p)}$ intersects $\lambda^{(p)}$ in exactly one point.

A consequence of the previous paragraph is: if $\mathcal D'$ is a
nonnegative relative periodic domain of $\Sigma-\widetilde A$,
then $\partial\mathcal D'$ does not pass through
$\overline{\xi}^{(p)}_i$'s. Moreover, $\partial\mathcal D'$ does
not pass through the rest of the curves in $\overline{A}^{(p)}$,
because
$\overline{\eta}^{(p)}_2,\dots,\overline{\eta}^{(p)}_{2g+1},\overline{\delta}^{(p)}$
are linearly independent in
$\homo_1(\overline{A}^{(p)},\partial\overline{A}^{(p)})$, and
$\alpha^{(p-1)}$ itself can not bound a relative periodic domain
in $\overline{A}^{(p)}$. Of course, here we can choose a point
near $\overline{\delta}^{(p)}\cap\overline{\lambda}^{(p)}$ as the
base point in $\overline{A}^{(p)}$.

We can go on with the above argument applied to
$A^{(p-1)},\overline A^{(p-1)},\dots,A^{(1)},B,\overline B$
inductively, to conclude that the local multiplicity of $\mathcal
D'$ is 0 in these subsurfaces.

We also wind
$\widetilde{\alpha}_{2g+2},\dots,\widetilde{\alpha}_{2g+1+r}$ in
$\widetilde{A}$, as in the proof of \fullref{support}.
Hence the new diagram
$(\Sigma^*,\mbox{\boldmath${\alpha}$}^{*\prime},\mbox{\boldmath${\beta}$}^*,w^*,z^*)$
after winding is weakly admissible. Moreover, the domains of the
holomorphic disks corresponding to the boundary map will restrict
to relative periodic domains in $\Sigma-\widetilde{A}$, so the
holomorphic disks for $\widehat{CFK}(Y^*,K^*,-(pg+h))$ are
supported inside $A^{(p)}\cup\widetilde A$. Hence they are in
one-to-one correspondence with the holomorphic disks for
$\widehat{CFK}(Y,K,-g)$. Then our desired result holds.
\end{proof}

\bibliographystyle{gtart}
\bibliography{link}

\begin{thebibliography}{}
\providecommand\bibmarginpar{\leavevmode\marginpar}
\def\urlstyle#1{{\tt #1}}

\bibitem{BZ}
\textbf{G Burde}, \textbf{H Zieschang}, \emph{Knots}, de Gruyter Studies in
  Mathematics 5, Walter de Gruyter \& Co., Berlin (2003) \xox{MR}{1959408}

\bibitem{Ef}
\textbf{E Eftekhary}, \emph{Longitude {F}loer homology and the {W}hitehead
  double}, Algebr. Geom. Topol. 5 (2005) 1389--1418 \xox{MR}{2171814}

\bibitem{G1}
\textbf{D Gabai}, \emph{The {M}urasugi sum is a natural geometric operation},
  from: ``Low-dimensional topology (San Francisco, Calif., 1981)'', Contemp.
  Math. 20, Amer. Math. Soc., Providence, RI (1983)  131--143 \xox{MR}{718138}

\bibitem{G2}
\textbf{D Gabai}, \emph{The {M}urasugi sum is a natural geometric operation.
  {II}}, from: ``Combinatorial methods in topology and algebraic geometry
  (Rochester, N.Y., 1982)'', Contemp. Math. 44, Amer. Math. Soc., Providence,
  RI (1985)  93--100 \xox{MR}{813105}

\bibitem{G3}
\textbf{D Gabai}, \emph{Foliations and the topology of {$3$}-manifolds. {III}},
  J. Differential Geom. 26 (1987) 479--536 \xox{MR}{910018}

\bibitem{He}
\textbf{M Hedden}, \href{http://dx.doi.org/10.2140/agt.2005.5.1197} {\emph{On
  knot {F}loer homology and cabling}}, Algebr. Geom. Topol. 5 (2005) 1197--1222
  \xox{MR}{2171808}

\bibitem{Ni}
\textbf{Y Ni}, \emph{A note on knot Floer homology of links}
  \xox{arXiv}{math.GT/0506208}

\bibitem{OSz9}
\textbf{P Ozsv{\'a}th}, \textbf{Z Szab{\'o}}, \emph{Heegaard diagrams and
  holomorphic disks}, from: ``Different faces of geometry'', Int. Math. Ser.
  (N. Y.), Kluwer/Plenum, New York (2004)  301--348 \xox{MR}{2102999}

\bibitem{OSz4}
\textbf{P Ozsv{\'a}th}, \textbf{Z Szab{\'o}},
  \href{http://dx.doi.org/10.2140/gt.2004.8.311} {\emph{Holomorphic disks and
  genus bounds}}, Geom. Topol. 8 (2004) 311--334 \xox{MR}{2023281}

\bibitem{OSz2}
\textbf{P Ozsv{\'a}th}, \textbf{Z Szab{\'o}},
  \href{http://dx.doi.org/10.1016/j.aim.2003.05.001} {\emph{Holomorphic disks
  and knot invariants}}, Adv. Math. 186 (2004) 58--116 \xox{MR}{2065507}

\bibitem{OSz1}
\textbf{P Ozsv{\'a}th}, \textbf{Z Szab{\'o}},
  \href{http://projecteuclid.org/getRecord?id=euclid.annm/1105737569}
  {\emph{Holomorphic disks and three-manifold invariants: properties and
  applications}}, Ann. of Math. $(2)$ 159 (2004) 1159--1245 \xox{MR}{2113020}

\bibitem{Osz0}
\textbf{P Ozsv{\'a}th}, \textbf{Z Szab{\'o}},
  \href{http://projecteuclid.org/getRecord?id=euclid.annm/1105737568}
  {\emph{Holomorphic disks and topological invariants for closed
  three-manifolds}}, Ann. of Math. $(2)$ 159 (2004) 1027--1158
  \xox{MR}{2113019}

\bibitem{OSz3}
\textbf{P Ozsv{\'a}th}, \textbf{Z Szab{\'o}},
  \href{http://projecteuclid.org/getRecord?id=euclid.dmj/1121448863}
  {\emph{Heegaard {F}loer homology and contact structures}}, Duke Math. J. 129
  (2005) 39--61 \xox{MR}{2153455}

\bibitem{Ra}
\textbf{J Rasmussen}, \emph{Floer homology and knot complements}, PhD thesis,
  Harvard (2003) \xox{arXiv}{math.GT/0306378}

\bibitem{Ro}
\textbf{D Rolfsen}, \emph{Knots and links}, Mathematics Lecture Series 7,
  Publish or Perish Inc., Houston, TX (1990) \xox{MR}{1277811}\ Corrected
  reprint of the 1976 original

\bibitem{St}
\textbf{J\,R Stallings}, \emph{Constructions of fibred knots and links}, from:
  ``Algebraic and geometric topology (Proc. Sympos. Pure Math., Stanford Univ.,
  Stanford, Calif., 1976), Part 2'', Proc. Sympos. Pure Math., XXXII, Amer.
  Math. Soc., Providence, R.I. (1978)  55--60 \xox{MR}{520522}

\end{thebibliography}

\end{document}